# LARGE DEVIATIONS FOR INFINITE DIMENSIONAL STOCHASTIC DYNAMICAL SYSTEMS

By Amarjit Budhiraja,[1] Paul Dupuis[2] and Vasileios Maroulas[1]

*University of North Carolina, Brown University and University of North Carolina*

The large deviations analysis of solutions to stochastic differential equations and related processes is often based on approximation. The construction and justification of the approximations can be onerous, especially in the case where the process state is infinite dimensional. In this paper we show how such approximations can be avoided for a variety of infinite dimensional models driven by some form of Brownian noise. The approach is based on a variational representation for functionals of Brownian motion. Proofs of large deviations properties are reduced to demonstrating basic qualitative properties (existence, uniqueness and tightness) of certain perturbations of the original process.

**1. Introduction.** Small noise large deviations theory for stochastic differential equations (SDE) has a long history. The finite dimensional setting, that is, where the SDE is driven by finitely many Brownian motions, was first studied by Freidlin and Wentzell [13]. In its basic form, one considers a $k$-dimensional SDE of the form

$$(1.1) \quad dX^\epsilon(t) = b(X^\epsilon(t))\,dt + \sqrt{\epsilon}a(X^\epsilon(t))\,dW(t), \qquad X^\varepsilon(0) = x^\varepsilon, t \in [0,T],$$

with coefficients $a$, $b$ satisfying suitable regularity properties and $W$ a finite dimensional standard Brownian motion. If $x^\varepsilon \to x^0$ as $\epsilon \to 0$, then $X^\epsilon \xrightarrow{\mathbb{P}} X^0$

Received March 2007; revised March 2007.

[1]Supported in part by the Army Research Office Grants W911NF-04-1-0230, W911NF-0-1-0080.

[2]Supported in part by the NSF Grants DMS-03-06070, DMS-04-04806 and DMS-07-06003 and the Army Research Office Grant W911NF-05-1-0289.

*AMS 2000 subject classifications.* Primary 60H15, 60F10; secondary 37L55.

*Key words and phrases.* Large deviations, Brownian sheet, Freidlin–Wentzell LDP, stochastic partial differential equations, stochastic evolution equations, small noise asymptotics, infinite dimensional Brownian motion.

in $\mathcal{C}([0,T]:\mathbb{R}^k)$, where $X^0$ solves the equation $\dot{x} = b(x)$ with initial data $x^0$. The Freidlin–Wentzell theory describes the path asymptotics, as $\varepsilon \to 0$, of probabilities of large deviations of the solution of the SDE from $X^0$—the law of large number dynamics. Since the original work of Freidlin–Wentzell, the finite dimensional problem has been extensively studied and many of the original assumptions made in [13] have been significantly relaxed (cf. [1, 8]).

Our interest in this work is with infinite dimensional models, that is, the setting where the driving Brownian motion $W$ is "infinite dimensional." In recent years there has been a lot of interest in large deviations analysis for such SDEs, and a partial list of references is [2, 4, 5, 6, 9, 12, 16, 19, 22, 25, 29, 30]. Our approach to the large deviation analysis, which is based on certain variational representations for infinite dimensional Brownian motions [3], is very different from that taken in these papers. The goal of the present work is to show how the variational representations can be easily applied to prove large deviation properties for diverse families of infinite dimensional models. Of course, the claim that a certain approach is easy to use may be viewed as subjective, and such a claim is only truly validated when other researchers find the approach convenient. In this regard, it is worth noting that the recent works [23, 24, 26, 27] have proved large deviation properties by applying the general large deviation principle (LDP) for Polish space valued measurable functionals of a Hilbert space valued Brownian motion established in [3] (see Section 6 for details).

As noted previously, one contribution of the present paper is to demonstrate in the context of an interesting example how easy it is to verify the main assumption for the LDP made in [3]. A second contribution is to show how the setup of [3], which considered SDEs driven by a Hilbert space valued Wiener processes, can be generalized to closely related settings, such as equations driven by a Brownian sheet. The chosen application is to a class of reaction-diffusion stochastic partial differential equations (SPDE) [see (5.1)], for which well-posedness has been studied in [20] and a small noise LDP established in [19]. The class includes, as a special case, the reaction-diffusion SPDEs considered in [25] (See Remark 3).

Our proof of the LDP proceeds by verification of the condition analogous to Assumption 4.3 of [3] (Assumption 3 in the current paper) appropriate to this formulation. The key ingredient in the verification of this assumption are the well-posedness and compactness for sequences of controlled versions of the original SPDE; see Theorems 10, 11 and 12. For comparison, the statements analogous to Theorems 10 and 11 in the finite dimensional setting (1.1) would say that, for any $\theta \in [0,1)$ and any $L^2$-bounded control $u$, [i.e., a predictable process satisfying $\int_0^T \|u(s)\|^2 \, ds \leq M$, a.s. for some $M \in (0,\infty)$], and any initial condition $x \in \mathbb{R}^k$, the equation

$$\begin{aligned}
dX_x^{\theta,u}(t) &= b(X_x^{\theta,u}(t))\,dt + \theta a(X_x^{\theta,u}(t))\,dW(t) + a(X_x^{\theta,u}(t))u(t)\,dt, \\
X_x^{\theta,u}(0) &= x
\end{aligned} \quad (1.2)$$

has a unique solution for $t \in [0, T]$. Also, the statement analogous to Theorem 12 in the finite dimensional setting would require that if $\theta(\varepsilon) \to \theta(0) = 0$, if a sequence of uniformly $L^2$-bounded controls $u^\varepsilon$ satisfies $u^\epsilon \to u$ in distribution (with the weak topology on the bounded $L^2$ ball), and if $x^{\theta(\epsilon)} \to x$ (all as $\varepsilon \to 0$), then $X^{\theta(\epsilon), u^\epsilon}_{x^{\theta(\epsilon)}} \to X^{0,u}_x$ in distribution.

As one may expect, the techniques and estimates used to prove such properties for the original (uncontrolled) stochastic model can be applied here as well, and indeed, proofs for the controlled SPDEs proceed in very much the same way as those of their uncontrolled counterparts. A side benefit of this pleasant situation is that one can often prove large deviation properties under mild conditions, and indeed, conditions that differ little from those needed for a basic qualitative analysis of the original equation. In the present setting, we are able to relax two of the main technical conditions used in [19], which are the uniform boundedness of the diffusion coefficient [i.e., the function $F$ in (5.1)] and the so-called "cone condition" imposed on the underlying domain (cf. [18], page 320). In place of these, we require only that the domain be a bounded open set and that the diffusion coefficients satisfy the standard linear growth condition. It is stated in Remark 3.2 of [19] that although unique solvability holds under the weaker linear growth condition, they are unable to derive the corresponding large deviation principle. The conditions imposed on $F$ and $\mathcal{O}$ in [19] enter in an important way in their proofs of the large deviation principle which is based on obtaining suitable exponential tail probability estimates for certain stochastic convolutions in Hölder norms. This relies on the application of a generalization of Garsia's theorem [14], which requires the restrictive conditions alluded to above. An important point is that these conditions are not needed for unique solvability of the SPDE.

In contrast, the weak convergence proof presented here does not require any exponential probability estimates and, hence, these assumptions are no longer needed. Indeed, suitable exponential continuity (in probability) and exponential tightness estimates are perhaps the hardest and most technical parts of the usual proofs based on discretization and approximation arguments. This becomes particularly hard in infinite dimensional settings where these estimates are needed with metrics on exotic function spaces (e.g., Hölder spaces, spaces of diffeomorphisms, etc.).

Standard approaches to small noise LDP for infinite dimensional SDE build on the ideas of [1]. The key ingredients to the proof are as follows. One first considers an approximating Gaussian model which is obtained from the original SDE by freezing the coefficients of the right-hand side according to a time discretization. Each such approximation is then further approximated by a finite dimensional system uniformly in the value of the frozen (state) variable. Next, one establishes an LDP for the finite dimensional system and argues that the LDP continues to hold as one approaches the infinite dimensional model. Finally, one needs to obtain suitable exponential continuity estimates in order to obtain the LDP for the original non-Gaussian model from that for the frozen Gaussian model. Exponential continuity (in probability) and exponential tightness estimates that are used to justify these approximations are often obtained under additional conditions on the model than those needed for well-posedness and compactness. In particular, as noted earlier, for the reaction diffusion systems considered here, these rely on exponential tail probability estimates in Hölder norms for certain stochastic convolutions which are only available for bounded integrands.

An alternative approach, based on nonlinear semigroup theory and infinite dimensional Hamilton–Jacobi (HJ) equations, has been developed in [10] (see also [11]). The method of proof involves showing that the value function of the limit control problem that is obtained by the law of large number analysis of certain controlled perturbations of the original stochastic model uniquely solves an appropriate infinite dimensional HJ equation in a suitable viscosity sense. In addition, one needs to establish exponential tightness by verifying a suitable exponential compact containment estimate. Although both these steps have been verified for a variety of models (cf. [11]), the proofs are quite technical and rely on a uniqueness theory for infinite dimensional nonlinear PDEs. The uniqueness requirement on the limit HJ equation is an extraneous artifact of the approach, and different stochastic models seem to require different methods for this, in general very hard, uniqueness problem. In contrast to the weak convergence approach, it requires an analysis of the model that goes significantly beyond the unique solvability of the SPDE. In addition, as discussed previously, the exponential tightness estimates are typically the most technical part of the large deviation analysis for infinite dimensional models, and are often only available under "sub-optimal" conditions when using standard techniques.

We now give an outline of the paper. Section 2 contains some background material on large deviations and infinite dimensional Brownian motions. We recall some basic definitions and the equivalence between a LDP and Laplace principle for a family of probability measures on some Polish space. We next recall some commonly used formulations for an infinite dimensional Brownian motion, such as an infinite sequence of i.i.d. standard real Brownian motions, a Hilbert space valued Brownian motion, a cylindrical Brownian motion and a space-time Brownian sheet. Relationships between these various formulations are noted as well. In Section 3 we present a variational representation for bounded nonnegative functionals of an infinite sequence of real Brownian motions. This variational representation, originally obtained in [3], is the starting point of our study. We also provide analogous representations for other formulations of infinite dimensional Brownian motions. Section 4 gives a general uniform large deviation result for Polish space valued functionals of an infinite dimensional Brownian motion. We provide sufficient conditions for the uniform LDP for each of the formulations of an infinite dimensional Brownian motion mentioned above. In Section 5 we introduce the small noise reaction-diffusion SPDE and use the general uniform LDP of Section 4 to establish a Freidlin–Wentzell LDP for such SPDEs in an appropriate Hölder space. Finally, Section 6 gives a brief overview of some other recent works that have used this variational approach to establish small noise LDP for infinite dimensional models. An Appendix collects proofs that are postponed for purposes of presentation.

Some notation and mathematical conventions used in this work are as follows. Infima over the empty set are taken to be $+\infty$. All Hilbert spaces in this work will be separable. The Borel sigma-field on a Polish space $\mathcal{S}$ will be denoted by $\mathcal{B}(\mathcal{S})$.

**2. Preliminaries.** In this section we present some standard definitions and results from the theory of large deviations and infinite dimensional Brownian motions.

*Large deviation principle and Laplace asymptotics.* Let $\{X^\epsilon, \epsilon > 0\} \equiv \{X^\epsilon\}$ be a family of random variables defined on a probability space $(\Omega, \mathcal{F}, \mathbb{P})$ and taking values in a Polish space (i.e., a complete separable metric space) $\mathcal{E}$. Denote the metric on $\mathcal{E}$ by $d(x, y)$ and expectation with respect to $\mathbb{P}$ by $\mathbb{E}$. The theory of large deviations is concerned with events $A$ for which probabilities $\mathbb{P}(X^\epsilon \in A)$ converge to zero exponentially fast as $\epsilon \to 0$. The exponential decay rate of such probabilities is typically expressed in terms of a "rate function" $I$ mapping $\mathcal{E}$ into $[0, \infty]$.

DEFINITION 1 (*Rate function*). A function $I: \mathcal{E} \to [0, \infty]$ is called a rate function on $\mathcal{E}$, if for each $M < \infty$ the level set $\{x \in \mathcal{E} : I(x) \leq M\}$ is a compact subset of $\mathcal{E}$. For $A \in \mathcal{B}(\mathcal{E})$, we define $I(A) \doteq \inf_{x \in A} I(x)$.

DEFINITION 2 (*Large deviation principle*). Let $I$ be a rate function on $\mathcal{E}$. The sequence $\{X^\epsilon\}$ is said to satisfy the large deviation principle on $\mathcal{E}$ with rate function I if the following two conditions hold:

1. *Large deviation upper bound.* For each closed subset F of $\mathcal{E}$,
$$\limsup_{\epsilon \to 0} \epsilon \log \mathbb{P}(X^\epsilon \in F) \leq -I(F).$$

2. *Large deviation lower bound.* For each open subset G of $\mathcal{E}$,
$$\liminf_{\epsilon \to 0} \epsilon \log \mathbb{P}(X^\epsilon \in G) \geq -I(G).$$

If a sequence of random variables satisfies the large deviation principle with some rate function, then the rate function is unique [8], Theorem 1.3.1. In many problems one is interested in obtaining exponential estimates on functions which are more general than indicator functions of closed or open sets. This leads to the study of the Laplace principle.

DEFINITION 3 (*Laplace principle*). Let $I$ be a rate function on $\mathcal{E}$. The sequence $\{X^\epsilon\}$ is said to satisfy the Laplace principle upper bound (resp. lower bound) on $\mathcal{E}$ with rate function $I$ if for all bounded continuous functions $h:\mathcal{E}\to\mathbb{R}$,

$$(2.1) \qquad \limsup_{\epsilon\to 0} \epsilon \log \mathbb{E}\left\{\exp\left[-\frac{1}{\epsilon}h(X^\epsilon)\right]\right\} \leq -\inf_{f\in\mathcal{E}}\{h(f)+I(f)\}$$

and, respectively,

$$(2.2) \qquad \liminf_{\epsilon\to 0} \epsilon \log \mathbb{E}\left\{\exp\left[-\frac{1}{\epsilon}h(X^\epsilon)\right]\right\} \geq -\inf_{f\in\mathcal{E}}\{h(f)+I(f)\}.$$

The Laplace principle is said to hold for $\{X^\epsilon\}$ with rate function $I$ if both the Laplace upper and lower bounds are satisfied for all bounded continuous functions $h$.

One of the main results of the theory of large deviations is the equivalence between the Laplace principle and the large deviation principle. For a proof we refer the reader to [8], Section 1.2.

THEOREM 1. *The family $\{X^\epsilon\}$ satisfies the Laplace principle upper (resp. lower) bound with a rate function $I$ on $\mathcal{E}$ if and only if $\{X^\epsilon\}$ satisfies the large deviation upper (resp. lower) bound for all closed sets (respectively open sets) with the rate function $I$.*

In view of this equivalence, the rest of this work will be concerned with the study of the Laplace principle. In fact, we will study a somewhat strengthened notion, namely, a *Uniform Laplace Principle*, as introduced below. The uniformity is critical in certain applications, such as the study of exit time and invariant measure asymptotics for small noise Markov processes [13].

Let $\mathcal{E}_0$ and $\mathcal{E}$ be Polish spaces. For each $\epsilon>0$ and $y\in\mathcal{E}_0$, let $X^{\epsilon,y}$ be $\mathcal{E}$ valued random variables given on the probability space $(\Omega,\mathcal{F},\mathbb{P})$.

DEFINITION 4. A family of rate functions $I_y$ on $\mathcal{E}$, parametrized by $y\in\mathcal{E}_0$, is said to have compact level sets on compacts if for all compact subsets $K$ of $\mathcal{E}_0$ and each $M<\infty$, $\Lambda_{M,K}\doteq\bigcup_{y\in K}\{x\in\mathcal{E}:I_y(x)\leq M\}$ is a compact subset of $\mathcal{E}$.

DEFINITION 5 (*Uniform Laplace Principle*). Let $I_y$ be a family of rate functions on $\mathcal{E}$ parameterized by $y$ in $\mathcal{E}_0$ and assume that this family has compact level sets on compacts. The family $\{X^{\epsilon,y}\}$ is said to satisfy the Laplace principle on $\mathcal{E}$ with rate function $I_y$, uniformly on compacts, if for all compact subsets $K$ of $\mathcal{E}_0$ and all bounded continuous functions $h$ mapping $\mathcal{E}$ into $\mathbb{R}$,

$$\lim_{\epsilon \to 0} \sup_{y \in K} \left| \epsilon \log \mathbb{E}_y \left\{ \exp\left[ -\frac{1}{\epsilon} h(X^{\epsilon,y}) \right] \right\} + \inf_{x \in \mathcal{E}} \{h(x) + I_y(x)\} \right| = 0.$$

We next summarize some well-known formulations for infinite dimensional Brownian motions and note some elementary relationships between them.

*Infinite dimensional Brownian motions.* An infinite dimensional Brownian motion arises in a natural fashion in the study of stochastic processes with a spatial parameter. We refer the reader to [7, 18, 28] for numerous examples in the physical sciences where an infinite dimensional Brownian motion is used to model the driving noise for some dynamical system. Depending on the application of interest, the infinite dimensional nature of the driving noise may be expressed in a variety of forms. Some examples include an infinite sequence of i.i.d. standard (1-dim) Brownian motions, a Hilbert space valued Brownian motion, a cylindrical Brownian motion and a space-time Brownian sheet. In what follows, we describe all of these models and explain how they are related to each other. We will be only concerned with processes defined over a fixed time horizon and thus fix a $T > 0$, and all filtrations and stochastic processes will be defined over the horizon $[0,T]$. Reference to $T$ will be omitted unless essential. Let $(\Omega, \mathcal{F}, \mathbb{P})$ be a probability space with an increasing family of right continuous $\mathbb{P}$-complete sigma fields $\{\mathcal{F}_t\}$. We refer to $(\Omega, \mathcal{F}, \mathbb{P}, \{\mathcal{F}_t\})$ as a filtered probability space.

Let $\{\beta_i\}_{i=1}^\infty$ be an infinite sequence of independent, standard, one dimensional, $\{\mathcal{F}_t\}$-Brownian motions given on this filtered probability space. We denote the product space of countably infinite copies of the real line by $\mathbb{R}^\infty$. Endowed with the topology of coordinate-wise convergence, $\mathbb{R}^\infty$ is a Polish space. Then $\beta = \{\beta_i\}_{i=1}^\infty$ is a random variable with values in the Polish space $\mathcal{C}([0,T]:\mathbb{R}^\infty)$ and represents the simplest model of an infinite dimensional Brownian motion.

Frequently in applications it is convenient to endow the state space of the driving noise, as in the finite dimensional theory, with an inner product structure. Let $(H, \langle \cdot, \cdot \rangle)$ be a real separable Hilbert space. Let $Q$ be a bounded, strictly positive, trace class operator on $H$.

DEFINITION 6. An $H$ valued stochastic process $\{W(t)\}$ defined on $(\Omega, \mathcal{F}, \mathbb{P}, \{\mathcal{F}_t\})$ is called a $Q$-Wiener process with respect to $\{\mathcal{F}_t\}$ if, for every nonzero $h \in H$,

$$\{\langle Qh, h\rangle^{-1/2} \langle W(t), h\rangle, \{\mathcal{F}_t\}\}$$

is a one-dimensional standard Wiener process.

It can be shown that if $W$ is an $H$ valued $Q$-Wiener process, then $\mathbb{P}[W \in \mathcal{C}([0,T]:H)] = 1$, where $\mathcal{C}([0,T]:H)$ is the space of continuous functions from the closed interval $[0,T]$ to the Hilbert space $H$. Let $\{e_i\}_{i=1}^\infty$ be a complete orthonormal system (CONS) for the Hilbert space $H$ such that $Qe_i = \lambda_i e_i$, where $\lambda_i$ is the strictly positive $i$th eigenvalue of $Q$ that corresponds to the eigenvector $e_i$. Since $Q$ is a trace class operator, $\sum_{i=1}^\infty \lambda_i < \infty$. Define $\tilde{\beta}_i(t) \doteq \langle W(t), e_i\rangle, t \geq 0, i \in \mathbb{N}$. It is easy to check that $\{\tilde{\beta}_i\}$ is a sequence of independent $\{\mathcal{F}_t\}$-Brownian motions with quadratic variation $\langle\langle \tilde{\beta}_i, \tilde{\beta}_j\rangle\rangle_t = \lambda_i \delta_{ij} t$, where $\delta_{ij} = 1$ if $i = j$ and 0 otherwise. Setting $\beta_i = \tilde{\beta}_i/\sqrt{\lambda_i}$, $\{\beta_i\}_{i=1}^\infty$ is a sequence of independent, standard, one dimensional, $\{\mathcal{F}_t\}$-Brownian motions. Thus, starting from a $Q$-Wiener process, one can produce an infinite collection of independent, standard Brownian motions in a straight forward manner. Conversely, given a collection of independent, standard Brownian motions $\{\beta_i\}_{i=1}^\infty$ and $(Q, \{e_i, \lambda_i\})$ as above, one can obtain a $Q$-Wiener process $W$ by setting

$$(2.3) \qquad W(t) \doteq \sum_{i=1}^\infty \sqrt{\lambda_i} \beta_i(t) e_i.$$

The right-hand side of (2.3) clearly converges in $L^2(\Omega)$ for each fixed $t$. Furthermore, one can check that the series also converges in $\mathcal{C}([0,T]:H)$ almost surely (see [7], Theorem 4.3). These observations lead to the following result.

PROPOSITION 1. *There exist measurable maps $f:\mathcal{C}([0,T]:\mathbb{R}^\infty) \longmapsto \mathcal{C}([0,T]:H)$ and $g:\mathcal{C}([0,T]:H) \longmapsto \mathcal{C}([0,T]:\mathbb{R}^\infty)$ such that $f(\beta) = W$ and $g(W) = \beta$ a.s.*

REMARK 1. Consider the Hilbert space $l_2 \doteq \{x \equiv (x_1, x_2, \ldots): x_i \in \mathbb{R}$ and $\sum x_i^2 < \infty\}$ with the inner product $\langle x, y\rangle \doteq \sum x_i y_i$. Let $\{\lambda_i\}_{i=1}^\infty$ be a sequence of strictly positive numbers such that $\sum \lambda_i < \infty$. Then the Hilbert space $\bar{l}_2 \doteq \{x \equiv (x_1, x_2, \ldots): x_i \in \mathbb{R}$ and $\sum \lambda_i x_i^2 < \infty\}$ with the inner product $\langle x, y\rangle_1 \doteq \sum \lambda_i x_i y_i$ contains $l_2$ and the embedding map is Hilbert–Schmidt. Furthermore, the infinite sequence of real Brownian motions $\beta$ takes values in $\bar{l}_2$ almost surely and can be regarded as a $\bar{l}_2$ valued $Q$-Wiener process with $\langle Qx, y\rangle_1 = \sum_{i=1}^\infty \lambda_i^2 x_i y_i$.

Equation (2.3) above can be interpreted as saying that the sequence $\{\lambda_i\}$ (or, equivalently, the trace class operator $Q$) injects a "coloring" to a white noise such that the resulting process has better regularity. In some models of interest, such coloring is obtained indirectly in terms of (state dependent) diffusion coefficients. It is natural in such situations to consider the driving noise as a "cylindrical Brownian motion" rather than a Hilbert space valued Brownian motion. Let $(H, \langle \cdot, \cdot \rangle)$ be a real separable Hilbert space and fix a filtered probability space as above.

DEFINITION 7. A family $\{B_t(h) \equiv B(t,h) : t \in [0,T], h \in H\}$ of real random variables is said to be an $\{\mathcal{F}_t\}$-cylindrical Brownian motion if:

1. For every $h \in H$ with $\|h\| = 1$, $\{B(t,h), \mathcal{F}_t\}$ is a standard Wiener process.
2. For every $t \geq 0, a_1, a_2 \in \mathbb{R}$ and $f_1, f_2 \in H$,

$$B(t, a_1 f_1 + a_2 f_2) = a_1 B(t, f_1) + a_2 B(t, f_2) \quad \text{a.s.}$$

Note that if $\{B_t(h) : t \geq 0, h \in H\}$ is a cylindrical Brownian motion and $\{e_i\}$ is a CONS in $H$, then setting $\beta_i(t) \doteq B(t, e_i)$, we see that $\{\beta_i\}$ is a sequence of independent, standard, real valued Brownian motions. Conversely, given a sequence $\{\beta_i\}_{i=1}^{\infty}$ of independent, standard Brownian motions on a filtered probability space,

$$(2.4) \qquad B_t(h) \doteq \sum_{i=1}^{\infty} \beta_i(t) \langle e_i, h \rangle$$

defines a cylindrical Brownian motion on $H$. For each $h \in H$, the series in (2.4) converges in $L^2(\Omega)$ and a.s. in $\mathcal{C}([0,T] : \mathbb{R})$.

PROPOSITION 2. *Let $B$ be a cylindrical Brownian motion as in Definition 7 and let $\beta$ be as constructed above. Then $\sigma\{B_s(h) : 0 \leq s \leq t, h \in H\} = \sigma\{\beta(s) : 0 \leq s \leq t\}$. In particular, if $X$ is a $\sigma\{B(s,h) : 0 \leq s \leq T, h \in H\}$ measurable random variable, then there exists a measurable map $g : \mathcal{C}([0,T] : \mathbb{R}^{\infty}) \longmapsto \mathbb{R}$ such that $g(\beta) = X$ a.s.*

In many physical dynamical systems with randomness, the driving noise is given as a space-time white noise process, also referred to as a Brownian sheet. In what follows we introduce this stochastic process and describe its relationship with the formulations considered above. Let $(\Omega, \mathcal{F}, \mathbb{P}, \{\mathcal{F}_t\})$ be a filtered probability space as before and fix a bounded open subset $\mathcal{O} \subseteq \mathbb{R}^d$.

DEFINITION 8. A Gaussian family of real valued random variables $\{B(t,x), (t,x) \in [0,T] \times \mathcal{O}\}$ on a filtered probability space is called a Brownian sheet if the following hold:

1. If $(t,x) \in [0,T] \times \mathcal{O}$, then $\mathbb{E}B(t,x) = 0$.
2. If $0 \leq s \leq t \leq T$ and $x \in \mathcal{O}$, then $B(t,x) - B(s,x)$ is independent of $\{\mathcal{F}_s\}$.
3. $\text{Cov}(B(t,x), B(s,y)) = \lambda(A_{t,x} \cap A_{s,y})$, where $\lambda$ is the Lebesgue measure on $[0,T] \times \mathcal{O}$ and $A_{t,x} \doteq \{(s,y) \in \mathbb{R}_+ \times \mathcal{O} : 0 \leq s \leq t \text{ and } y_j \leq x_j, j = 1, \ldots, d\}$.
4. The map $(t,u) \mapsto B(t,u)$ from $[0,T] \times \mathcal{O}$ to $\mathbb{R}$ is continuous a.s.

To introduce stochastic integrals with respect to a Brownian sheet, we need the following definitions.

DEFINITION 9 (*Elementary and simple functions*). A function $f : \mathcal{O} \times [0,T] \times \Omega \to \mathbb{R}$ is elementary if there exist $a, b \in [0,T], a \leq b$, a bounded $\{\mathcal{F}_a\}$-measurable random variable $X$ and $A \in \mathcal{B}(\mathcal{O})$ such that
$$f(x,s,\omega) = X(\omega)\mathbf{1}_{(a,b]}(s)\mathbf{1}_A(x).$$
A finite sum of elementary functions is referred to as a simple function. We denote by $\overline{\mathcal{S}}$ the class of all simple functions.

DEFINITION 10 (*Predictable $\sigma$-field*). The predictable $\sigma$-field $\mathcal{P}$ on $\Omega \times [0,T] \times \mathcal{O}$ is the $\sigma$-field generated by $\overline{\mathcal{S}}$. A function $f : \Omega \times [0,T] \times \mathcal{O} \to \mathbb{R}$ is called a predictable process if it is $\mathcal{P}$-measurable.

Let $\mathcal{P}_2$ be the class of all predictable processes $f$ such that $\int_{[0,T] \times \mathcal{O}} f^2(s,x) \, ds \, dx$ is finite a.s. Also, let $\mathcal{L}_2$ be the subset of those processes that satisfy $\int_{[0,T] \times \mathcal{O}} \mathbb{E} f^2(s,x) \, ds \, dx < \infty$. For all $f \in \mathcal{P}_2$, the stochastic integral $M_t(f) \doteq \int_{[0,t] \times \mathcal{O}} f(s,u) B(ds\,du)$, $t \in [0,T]$ is well defined as in Chapter 2 of [28]. Furthermore, for all $f \in \mathcal{P}_2$, $\{M_t(f)\}_{0 \leq t \leq T}$ is a continuous $\{\mathcal{F}_t\}$-local martingale which is in fact a square integrable martingale if $f \in \mathcal{L}_2$. The quadratic variation of this local martingale is given as $\langle\!\langle M(f), M(f) \rangle\!\rangle_t \doteq \int_{[0,t] \times \mathcal{O}} f^2(s,x) \, ds \, dx$. More properties of the stochastic integral can be found in [28].

Let $\{\phi_i\}_{i=1}^\infty$ be a CONS in $L^2(\mathcal{O})$. Then it is easy to verify that $\beta \equiv \{\beta_i\}_{i=1}^\infty$ defined as $\beta_i(t) \doteq \int_{[0,t] \times \mathcal{O}} \phi_i(x) B(ds\,dx), i \geq 1, t \in [0,T]$ is a sequence of independent, standard, real Brownian motions. Also for $(t,x) \in [0,T] \times \mathcal{O}$,

$$(2.5) \qquad B(t,x) = \sum_{i=1}^\infty \beta_i(t) \int_{\mathcal{O}} \phi_i(y) \mathbf{1}_{(-\infty,x]}(y) \, dy,$$

where $(-\infty, x] = \{y : y_i \leq x_i \text{ for all } i = 1, \ldots, d\}$ and the series in (2.5) converges in $L^2(\Omega)$ for each $(t,x)$. From these considerations it follows that

$$(2.6) \qquad \sigma\{B(t,x), t \in [0,T], x \in \mathcal{O}\} = \sigma\{\beta_i(t), i \geq 1, t \in [0,T]\}.$$

As a consequence of (2.6), we have the following result.

PROPOSITION 3. *There exists a measurable map $g : \mathcal{C}([0,T] : \mathbb{R}^\infty) \to \mathcal{C}([0,T] \times \mathcal{O} : \mathbb{R})$ such that $B = g(\beta)$ a.s., where $\beta$ is as defined above (2.5).*

**3. Variational representations.** The large deviation results established in this work critically use certain variational representations for infinite dimensional Brownian motions. Let $(\Omega, \mathcal{F}, \mathbb{P}, \{\mathcal{F}_t\})$ be as before and let $\beta = \{\beta_i\}$ be a sequence of independent real standard Brownian motions. Recall that $\beta$ is a $\mathcal{C}([0,T]:\mathbb{R}^\infty)$ valued random variable. We call a function $f:[0,T] \times \Omega \to \mathbb{R}$ *elementary* if there exist $a, b \in [0,T], a \leq b$, and a bounded $\{\mathcal{F}_a\}$-measurable random variable $X$ such that $f(s,\omega) = X(\omega)\mathbf{1}_{(a,b]}(s)$. A finite sum of elementary functions is referred to as a simple function. We denote by $\overline{\mathcal{S}}$ the class of all simple functions. The predictable $\sigma$-field $\mathcal{P}$ on $\Omega \times [0,T]$ is the $\sigma$-field generated by $\overline{\mathcal{S}}$. For a Hilbert space $(H, \langle \cdot, \cdot \rangle)$, a function $f: \Omega \times [0,T] \to H$ is called an $H$ valued predictable process if it is $\mathcal{P}$-measurable. Let $\mathcal{P}_2(H)$ be the family of all $H$ valued predictable processes for which $\int_0^T \|\phi(s)\|^2\, ds < \infty$ a.s., where $\|\cdot\|$ is the norm in the Hilbert space $H$. Note that in the case $H = l_2$, $u \in \mathcal{P}_2(H) = \mathcal{P}_2(l_2)$ can be written as $u = \{u_i\}_{i=1}^\infty$, where $u_i \in \mathcal{P}_2(\mathbb{R})$ and $\sum_{i=1}^\infty \int_0^T |u_i(s)|^2\, ds < \infty$ a.s.

THEOREM 2. *Let $\|\cdot\|$ denote the norm in the Hilbert space $l_2$ and let $f$ be a bounded, Borel measurable function mapping $\mathcal{C}([0,T]:\mathbb{R}^\infty)$ into $\mathbb{R}$. Then,*

$$-\log \mathbb{E}(\exp\{-f(\beta)\}) = \inf_{u \in \mathcal{P}_2(l_2)} \mathbb{E}\left( \tfrac{1}{2} \int_0^T \|u(s)\|^2\, ds + f\left(\beta + \int_0^\cdot u(s)\, ds\right)\right).$$

The representation established in [3] is stated in a different form but is equivalent to Theorem 2. Let $(H, \langle \cdot, \cdot \rangle)$ be a Hilbert space and let $W$ be an $H$ valued $Q$-Wiener process, where $Q$ is a bounded, strictly positive, trace class operator on the Hilbert space $H$. Let $H_0 = Q^{1/2}H$, then $H_0$ is a Hilbert space with the inner product $\langle h, k \rangle_0 \doteq \langle Q^{-1/2}h, Q^{-1/2}h \rangle, h, k \in H_0$. Also the embedding map $i: H_0 \mapsto H$ is a Hilbert–Schmidt operator and $ii^* = Q$. Let $\|\cdot\|_0$ denote the norm in the Hilbert space $H_0$. The following theorem is proved in [3]. Theorem 2 follows from Theorem 3 and Remark 1.

THEOREM 3. *Let $f$ be a bounded, Borel measurable function mapping $\mathcal{C}([0,T]:H)$ into $\mathbb{R}$. Then*

$$-\log \mathbb{E}(\exp\{-f(W)\}) = \inf_{u \in \mathcal{P}_2(H_0)} \mathbb{E}\left( \tfrac{1}{2} \int_0^T \|u(s)\|_0^2\, ds + f\left(W + \int_0^\cdot u(s)\, ds\right)\right).$$

We finally note the following representation theorem for a Brownian sheet which follows from Theorem 2, Proposition 3 and an application of Girsanov's theorem.

THEOREM 4. *Let $f : \mathcal{C}([0,T] \times \mathcal{O} : \mathbb{R}) \to \mathbb{R}$ be a bounded measurable map. Let $B$ be a Brownian sheet as in Definition* 8. *Then*

$$-\log \mathbb{E}(\exp\{-f(B)\}) = \inf_{u \in \mathcal{P}_2} \mathbb{E}\bigg(\frac{1}{2} \int_0^T \int_{\mathcal{O}} u^2(s,r) \, dr \, ds + f(B^u)\bigg),$$

*where* $B^u(t,x) = B(t,x) + \int_0^t \int_{(-\infty,x] \cap \mathcal{O}} u(s,y) \, dy \, ds$.

**4. Large deviations for functionals of infinite dimensional Brownian motions.** In this section we give sufficient conditions for the uniform Laplace principle for functionals of an infinite dimensional Brownian motion. The uniformity is with respect to a parameter $x$ (typically an initial condition), which takes values in some compact subset of a Polish space $\mathcal{E}_0$. The analogous nonuniform result was established in [3]. The proof for the uniform case uses only minor modifications, but for the sake of completeness we include the details in the Appendix.

We begin by considering the case of a Hilbert space valued Wiener process and then use this case to deduce analogous Laplace principle results for functionals of a cylindrical Brownian motion and a Brownian sheet. Let $(\Omega, \mathcal{F}, \mathbb{P}, \{\mathcal{F}_t\})$, $(H, \langle \cdot, \cdot \rangle)$, $Q$ be as in Section 2 and let $W$ be an $H$ valued Wiener process with trace class covariance $Q$ given on this filtered probability space (see Definition 6). Let $\mathcal{E}$ be a Polish space and for each $\epsilon > 0$, let $\mathcal{G}^\epsilon : \mathcal{E}_0 \times \mathcal{C}([0,T] : H) \to \mathcal{E}$ be a measurable map. We next present a set of sufficient conditions for a uniform large deviation principle to hold for the family $\{X^{\epsilon,x} \doteq \mathcal{G}^\epsilon(x, \sqrt{\epsilon}W)\}$ as $\epsilon \to 0$. Let $H_0$ be as introduced above Theorem 3 and define for $N \in \mathbb{N}$

$$(4.1) \qquad S^N(H_0) \doteq \bigg\{u \in L^2([0,T] : H_0) : \int_0^T \|u(s)\|_0^2 \, ds \leq N\bigg\},$$

$$(4.2) \qquad \mathcal{P}_2^N(H_0) \doteq \{u \in \mathcal{P}_2(H_0) : u(\omega) \in S^N(H_0), \mathbb{P}\text{-a.s.}\}.$$

It is easy to check that $S^N(H_0)$ is a compact metric space under the metric

$$d_1(x,y) = \sum_{i=1}^\infty \frac{1}{2^i} \bigg|\int_0^T \langle x(s) - y(s), e_i(s) \rangle_0 \, ds\bigg|.$$

Henceforth, wherever we refer to $S^N(H_0)$, we will consider it endowed with the topology obtained from the metric $d_1$ and refer to this as the weak topology on $S^N(H_0)$.

ASSUMPTION 1. *There exists a measurable map* $\mathcal{G}^0 : \mathcal{E}_0 \times \mathcal{C}([0,T] : H) \to \mathcal{E}$ *such that the following hold:*

1. *For every $M < \infty$ and compact set $K \subseteq \mathcal{E}_0$, the set*

$$\Gamma_{M,K} \doteq \bigg\{\mathcal{G}^0\bigg(x, \int_0^\cdot u(s) \, ds\bigg) : u \in S^M(H_0), x \in K\bigg\}$$

is a compact subset of $\mathcal{E}$.

2. Consider $M < \infty$ and families $\{u^\epsilon\} \subset \mathcal{P}_2^M(H_0)$ and $\{x^\varepsilon\} \subset \mathcal{E}_0$ such that $u^\epsilon$ converges in distribution [as $S^M(H_0)$ valued random elements] to $u$ and $x^\varepsilon \to x$ as $\epsilon \to 0$. Then

$$\mathcal{G}^\epsilon\left(x^\varepsilon, \sqrt{\epsilon}W(\cdot) + \int_0^\cdot u^\epsilon(s)\,ds\right) \to \mathcal{G}^0\left(x, \int_0^\cdot u(s)\,ds\right)$$

in distribution as $\epsilon \to 0$.

THEOREM 5. *Let $X^{\epsilon,x} = \mathcal{G}^\epsilon(x, \sqrt{\epsilon}W)$ and suppose that Assumption 1 holds. For $x \in \mathcal{E}_0$ and $f \in \mathcal{E}$, let*

$$(4.3) \qquad I_x(f) \doteq \inf_{\{u \in L^2([0,T]:H_0): f = \mathcal{G}^0(x,\int_0^\cdot u(s)ds)\}} \left\{\tfrac{1}{2}\int_0^T \|u(s)\|_0^2\,ds\right\}.$$

*Suppose that for all $f \in \mathcal{E}$, $x \mapsto I_x(f)$ is a lower semi-continuous (l.s.c.) map from $\mathcal{E}_0$ to $[0,\infty]$. Then, for all $x \in \mathcal{E}_0$, $f \mapsto I_x(f)$ is a rate function on $\mathcal{E}$ and the family $\{I_x(\cdot), x \in \mathcal{E}_0\}$ of rate functions has compact level sets on compacts. Furthermore, the family $\{X^{\epsilon,x}\}$ satisfies the Laplace principle on $\mathcal{E}$, with rate function $I_x$, uniformly on compact subsets of $\mathcal{E}_0$.*

As noted earlier, an analogous non-uniform result was established in [3]. We remark that there is a slight change in notation from [3]. Denoting the map $\mathcal{G}^\varepsilon$ introduced in [3] by $\overline{\mathcal{G}}^\varepsilon$, the correspondence with the $\mathcal{G}^\varepsilon$ introduced in this section is given as $\overline{\mathcal{G}}^\varepsilon(x,f) = \mathcal{G}^\varepsilon(x, f\sqrt{\varepsilon})$ for $x \in \mathcal{E}_0$ and $f \in \mathcal{C}([0,T]:H)$.

Next, let $\beta \equiv \{\beta_i\}$ be a sequence of independent standard real Brownian motions on $(\Omega, \mathcal{F}, \mathbb{P}, \{\mathcal{F}_t\})$. Recall that $\beta$ is a $(\mathcal{C}([0,T]:\mathbb{R}^\infty), \mathcal{B}(\mathcal{C}([0,T]:\mathbb{R}^\infty))) \equiv (S, \mathcal{S})$ valued random variable. For each $\varepsilon > 0$, let $\mathcal{G}^\varepsilon : \mathcal{E}_0 \times S \to \mathcal{E}$ be a measurable map and define

$$(4.4) \qquad X^{\varepsilon,x} \doteq \mathcal{G}^\varepsilon(x, \sqrt{\epsilon}\beta).$$

We now consider the Laplace principle for the family $\{X^{\varepsilon,x}\}$ and introduce the analog of Assumption 1 for this setting. In the assumption, $S^M(l_2)$ and $\mathcal{P}_2^M(l_2)$ are defined as in (4.1) and (4.2), with $H_0$ there replaced by the Hilbert space $l_2$.

ASSUMPTION 2. There exists a measurable map $\mathcal{G}^0 : \mathcal{E}_0 \times S \to \mathcal{E}$ such that the following hold:

1. For every $M < \infty$ and compact set $K \subseteq \mathcal{E}_0$, the set

$$\Gamma_{M,K} \doteq \left\{\mathcal{G}^0\left(x, \int_0^\cdot u(s)\,ds\right) u \in S^M(l_2), x \in K\right\}$$

is a compact subset of $\mathcal{E}$.

2. Consider $M < \infty$ and families $\{u^\epsilon\} \subset \mathcal{P}_2^M(l_2)$ and $\{x^\varepsilon\} \subset \mathcal{E}_0$ such that $u^\epsilon$ converges in distribution [as $S^M(l_2)$ valued random elements] to u and $x^\varepsilon \to x$ as $\epsilon \to 0$. Then

$$\mathcal{G}^\epsilon\left(x^\varepsilon, \sqrt{\epsilon}\beta + \int_0^\cdot u^\epsilon(s)\,ds\right) \to \mathcal{G}^0\left(x, \int_0^\cdot u(s)\,ds\right),$$

as $\varepsilon \to 0$ in distribution.

The proof of the following, which uses a straightforward reduction to Theorem 5, is given in the Appendix.

THEOREM 6. *Let $X^{\epsilon,x}$ be as in (4.4) and suppose that Assumption 2 holds. For $x \in \mathcal{E}_0$ and $f \in \mathcal{E}$, let*

$$(4.5) \qquad I_x(f) \doteq \inf_{\{u \in L^2([0,T]:l_2): f = \mathcal{G}^0(x, \int_0^\cdot u(s)ds)\}} \left\{\tfrac{1}{2}\int_0^T \|u(s)\|_{l_2}^2\,ds\right\}.$$

*Suppose that for all $f \in \mathcal{E}$, $x \mapsto I_x(f)$ is a l.s.c. map from $\mathcal{E}_0$ to $[0,\infty]$. Then, for all $x \in \mathcal{E}_0$, $f \mapsto I_x(f)$ is a rate function on $\mathcal{E}$ and the family $\{I_x(\cdot), x \in \mathcal{E}_0\}$ of rate functions has compact level sets on compacts. Furthermore, the family $\{X^{\epsilon,x}\}$ satisfies the Laplace principle on $\mathcal{E}$, with rate function $I_x$, uniformly on compact subsets of $\mathcal{E}_0$.*

Finally, to close this section, we consider the Laplace principle for functionals of a Brownian sheet. Let $B$ be a Brownian sheet as in Definition 8. Let $\mathcal{G}^\varepsilon : \mathcal{E}_0 \times \mathcal{C}([0,T] \times \mathcal{O}:\mathbb{R}) \to \mathcal{E}$, $\varepsilon > 0$ be a family of measurable maps. Define $X^{\varepsilon,x} \doteq \mathcal{G}^\varepsilon(x, \sqrt{\epsilon}B)$. We now provide sufficient conditions for the Laplace principle to hold for the family $\{X^{\varepsilon,x}\}$.

Analogous to classes defined in (4.1) and (4.2), we introduce

$$(4.6) \quad \begin{aligned} S^N &\doteq \left\{\phi \in L^2([0,T] \times \mathcal{O}) : \int_{[0,T] \times \mathcal{O}} \phi^2(s,r)\,ds\,dr \leq N\right\}, \\ \mathcal{P}_2^N &\doteq \{u \in \mathcal{P}_2 : u(\omega) \in S^N, \mathbb{P}\text{-a.s.}\}. \end{aligned}$$

Once more, $S^N$ is endowed with the weak topology on $L^2([0,T] \times \mathcal{O})$, under which it is a compact metric space. For $u \in L^2([0,T] \times \mathcal{O})$, define $\text{Int}(u) \in \mathcal{C}([0,T] \times \mathcal{O}:\mathbb{R})$ by

$$(4.7) \qquad \text{Int}(u)(t,x) \doteq \int_{[0,t] \times (\mathcal{O} \cap (-\infty,x])} u(s,y)\,ds\,dy,$$

where, as before, $(-\infty,x] = \{y : y_i \leq x_i \text{ for all } i=1,\ldots,d\}$.

ASSUMPTION 3. There exists a measurable map $\mathcal{G}^0 : \mathcal{E}_0 \times \mathcal{C}([0,T] \times \mathcal{O}:\mathbb{R}) \to \mathcal{E}$ such that the following hold:

1. For every $M < \infty$ and compact set $K \subseteq \mathcal{E}_0$, the set

$$\Gamma_{M,K} \doteq \{\mathcal{G}^0(x, \text{Int}(u)) : u \in S^M, x \in K\}$$

is a compact subset of $\mathcal{E}$, where $\text{Int}(u)$ is as defined in (4.7).
2. Consider $M < \infty$ and families $\{u^\varepsilon\} \subset \mathcal{P}_2^M$ and $\{x^\varepsilon\} \subset \mathcal{E}_0$ such that $u^\epsilon$ converges in distribution (as $S^M$ valued random elements) to $u$ and $x^\varepsilon \to x$ as $\epsilon \to 0$. Then

$$\mathcal{G}^\epsilon(x^\epsilon, \sqrt{\epsilon}B + \text{Int}(u^\epsilon)) \to \mathcal{G}^0(x, \text{Int}(u)),$$

in distribution as $\epsilon \to 0$.

For $f \in \mathcal{E}$ and $x \in \mathcal{E}_0$, define

$$(4.8) \quad I_x(f) = \inf_{\{u \in L^2([0,T] \times \mathcal{O}) : f = \mathcal{G}^0(x, \text{Int}(u))\}} \left\{ \tfrac{1}{2} \int_{[0,T] \times \mathcal{O}} u^2(s,r)\, dr\, ds \right\}.$$

THEOREM 7. *Let $\mathcal{G}^0 : \mathcal{E}_0 \times \mathcal{C}([0,T] \times \mathcal{O} : \mathbb{R}) \to \mathcal{E}$ be a measurable map satisfying Assumption 3. Suppose that for all $f \in \mathcal{E}, x \mapsto I_x(f)$ is a l.s.c. map from $\mathcal{E}_0$ to $[0, \infty]$. Then for every $x \in \mathcal{E}_0$, $I_x : \mathcal{E} \to [0, \infty]$, defined by (4.8), is a rate function on $\mathcal{E}$ and the family $\{I_x, x \in \mathcal{E}_0\}$ of rate functions has compact level sets on compacts. Furthermore, the family $\{X^{\epsilon,x}\}$ satisfies the Laplace principle on $\mathcal{E}$ with rate function $I_x$, uniformly for $x$ in compact subsets of $\mathcal{E}_0$.*

The proof of this theorem can be found in the Appendix.

## 5. Stochastic reaction-diffusion systems.

5.1. *The large deviation theorem.* In this section we will use results from Section 4, and in particular, Theorem 7, to study the small noise large deviations principle for a class of stochastic partial differential equations (SPDE) that has been considered in [20]. The class includes, as a special case, the reaction-diffusion SPDEs considered in [25] (see Remark 3). The main result of the section is Theorem 9, which establishes the uniform Freidlin–Wentzell LDP for such SPDEs.

As discussed previously, the weak convergence method bypasses the various discretizations, approximations and exponential probability estimates that are commonly used in standard approaches to the problem. Instead, one needs to only prove various qualitative properties (compactness, convergence, etc.) for sequences of controlled versions of the SPDE model. As one might expect, the techniques and estimates used to prove these properties for the original SPDE can be applied here as well, and indeed, proofs for the controlled SPDEs proceed in very much the same way as those of their uncontrolled counterparts.

Let $(\Omega, \mathcal{F}, \mathbb{P})$ be a probability space with an increasing family of right-continuous, $\mathbb{P}$-complete $\sigma$-fields $\{\mathcal{F}_t\}_{0 \leq t \leq T}$. Let $\mathcal{O} \subseteq \mathbb{R}^d$ be a bounded open set and $\{B(t,x) : (t,x) \in \mathbb{R}_+ \times \mathcal{O}\}$ be a Brownian sheet given on this filtered probability space. Consider the SPDE

$$
\begin{aligned}
dX(t,r) = {}& (L(t)X(t,r) + R(t,r,X(t,r)))\,dr\,dt \\
& + \sqrt{\epsilon} F(t,r,X(t,r)) B(dr\,dt)
\end{aligned}
\tag{5.1}
$$

with initial condition

$$X(0,r) = \xi(r).$$

Here $F$ and $R$ are measurable maps from $[0,T] \times \mathcal{O} \times \mathbb{R}$ to $\mathbb{R}$ and $\epsilon \in (0,\infty)$. Also, $\{L(t) : t \geq 0\}$ is a family of linear, closed, densely defined operators on $\mathcal{C}(\mathcal{O})$ that generates a two parameter strongly continuous semigroup $\{U(t,s) : 0 \leq s \leq t\}$ on $\mathcal{C}(\mathcal{O})$, with kernel function $G(t,s,r,q), 0 \leq s < t, r,q \in \mathcal{O}$. Thus, for $f \in \mathcal{C}(\mathcal{O})$, $U(t,t)f = f$, $t \in [0,T]$ and

$$(U(t,s)f)(r) = \int_\mathcal{O} G(t,s,r,q) f(q)\,dq, \qquad r \in \mathcal{O}, 0 \leq s < t \leq T.$$

For notational convenience, we write $f(r) = \int_\mathcal{O} G(0,0,r,q) f(q)\,dq$ for $f \in \mathcal{C}(\mathcal{O})$.

By a solution of the SPDE (5.1), we mean the following:

DEFINITION 11. A random field $X \equiv \{X(t,r) : t \in [0,T], r \in \mathcal{O}\}$ is called a mild solution of the stochastic partial differential equation (5.1) with initial condition $\xi$ if $(t,r) \mapsto X(t,r)$ is continuous a.s., $X(t,r)$ is $\{\mathcal{F}_t\}$-measurable for any $t \in [0,T], r \in \mathcal{O}$, and if

$$
\begin{aligned}
X(t,r) = {}& \int_\mathcal{O} G(t,0,r,q) \xi(q)\,dq \\
& + \int_0^t \int_\mathcal{O} G(t,s,r,q) R(s,q,X(s,q))\,dq\,ds \\
& + \sqrt{\epsilon} \int_0^t \int_\mathcal{O} G(t,s,r,q) F(s,q,X(s,q)) B(dq\,ds) \qquad \text{a.s.}
\end{aligned}
\tag{5.2}
$$

Implicit in Definition 11 is the requirement that the integrals in (5.2) are well defined. We will shortly introduce conditions on $G, F$ and $R$ that ensure that for a continuous adapted random field $X$, all the integrals in (5.2) are meaningful. As a convention, we take $G(t,s,r,q)$ to be zero when $0 \leq t \leq s \leq T, r,q \in \mathcal{O}$.

For $u \in \mathcal{P}_2^N$ [which was defined in (4.6)], the controlled analogue of (5.2) is

$$
\begin{aligned}
Y(t,r) &= \int_{\mathcal{O}} G(t,0,r,q)\xi(q)\,dq \\
&\quad + \int_0^t \int_{\mathcal{O}} G(t,s,r,q)R(s,q,Y(s,q))\,dq\,ds \\
&\quad + \sqrt{\epsilon}\int_0^t \int_{\mathcal{O}} G(t,s,r,q)F(s,q,Y(s,q))B(dq\,ds) \\
&\quad + \int_0^t \int_{\mathcal{O}} G(t,s,r,q)F(s,q,Y(s,q))u(s,q)\,dq\,ds.
\end{aligned}
\tag{5.3}
$$

As discussed previously, the main work in proving an LDP for (5.2) will be to prove qualitative properties (existence and uniqueness, tightness properties, and stability under perturbations) for solutions to (5.3). We begin by discussing known qualitative theory for (5.2).

For $\alpha > 0$, let $\mathbb{B}_\alpha = \{\psi \in \mathcal{C}(\mathcal{O}) : \|\psi\|_\alpha < \infty\}$ be the Banach space with norm

$$\|\psi\|_\alpha = \|\psi\|_0 + \sup_{r,q \in \mathcal{O}} \frac{|\psi(r) - \psi(q)|}{|r-q|^\alpha},$$

where $\|\psi\|_0 = \sup_{r \in \mathcal{O}} |\psi(r)|$. The Banach space $\mathbb{B}_\alpha([0,T] \times \mathcal{O})$ is defined similarly and for notational convenience, we denote this space by $\mathbb{B}_\alpha^T$. For $\alpha = 0$, the space $\mathbb{B}_0^T$ is the space of all continuous maps from $[0,T] \times \mathcal{O}$ to $\mathbb{R}$ endowed with the sup-norm. The following will be a standing assumption for this section. In the assumption, $\bar{\alpha}$ is a fixed constant, and the large deviation principle will be proved in the topology of $\mathcal{C}([0,T] : \mathbb{B}_\alpha)$, for any fixed $\alpha \in (0, \bar{\alpha})$.

ASSUMPTION 4. The following two conditions hold:

1. There exist constants $K(T) < \infty$ and $\gamma \in (d, \infty)$ such that:
   (a) for all $t, s \in [0, T], r \in \mathcal{O}$,
   $$\int_{\mathcal{O}} |G(t,s,r,q)|\,dq \leq K(T), \tag{5.4}$$
   (b) for all $0 \leq s < t \leq T$ and $r, q \in \mathcal{O}$,
   $$|G(t,s,r,q)| \leq K(T)(t-s)^{-d/\gamma}, \tag{5.5}$$
   (c) if $\bar{\alpha} = \frac{\gamma - d}{2\gamma}$, then for any $\alpha \in (0, \bar{\alpha})$ and for all $0 \leq s < t_1 \leq t_2 \leq T, r_1, r_2, q \in \mathcal{O}$,
   $$
   \begin{aligned}
   &|G(t_1,s,r_1,q) - G(t_2,s,r_2,q)| \\
   &\leq K(T)[(t_2-t_1)^{1-d/\gamma}(t_1-s)^{-1} + |r_1-r_2|^{2\alpha}(t_1-s)^{-(d+2\alpha)/\gamma}],
   \end{aligned}
   \tag{5.6}
   $$

(d) for all $x, y \in \mathbb{R}$, $r \in \mathcal{O}$ and $0 \leq t \leq T$,

(5.7) $$|R(t,r,x) - R(t,r,y)| + |F(t,r,x) - F(t,r,y)| \leq K(T)|x-y|$$

and

(5.8) $$|R(t,r,x)| + |F(t,r,x)| \leq K(T)(1+|x|).$$

2. For any $\alpha \in (0, \bar{\alpha})$ and $\xi \in \mathbb{B}_\alpha$, the trajectory $t \mapsto \int_\mathcal{O} G(t, 0, \cdot, q)\xi(q)\,dq$ belongs to $\mathcal{C}([0,T]:\mathbb{B}_\alpha)$ and the map

$$\mathbb{B}_\alpha \ni \xi \longmapsto \left\{ t \mapsto \int_\mathcal{O} G(t, 0, \cdot, q)\xi(q)\,dq \right\} \in \mathcal{C}([0,T]:\mathbb{B}_\alpha)$$

is a continuous map.

For future reference we recall that $\bar{\alpha} = \frac{\gamma - d}{2\gamma}$ and note that $\bar{\alpha} \in (0, 1/2)$.

REMARK 2.   1. We refer the reader to [19] for examples of families $\{L(t)\}_{t \geq 0}$ that satisfy this assumption.

2. Using (5.4) and (5.5), it follows that, for any $0 \leq s < t \leq T$ and $r \in \mathcal{O}$,

(5.9) $$\int_\mathcal{O} |G(t, s, r, q)|^2\,dq \leq K^2(T)(t-s)^{-d/\gamma}.$$

This, in particular, ensures that the stochastic integral in (5.2) is well defined.

3. Lemma 4.1(ii) of [19] shows that, under Assumption 4, for any $\alpha < \bar{\alpha}$ there exists a constant $\tilde{K}(\alpha)$ such that, for all $0 \leq t_1 \leq t_2 \leq T$ and all $r_1, r_2 \in \mathcal{O}$,

(5.10) $$\int_0^T \int_\mathcal{O} |G(t_1, s, r_1, q) - G(t_2, s, r_2, q)|^2\,dq\,ds$$
$$\leq \tilde{K}(\alpha)\rho((t_1, r_1), (t_2, r_2))^{2\alpha},$$

where $\rho$ is the Euclidean distance in $[0,T] \times \mathcal{O} \subset \mathbb{R}^{d+1}$. This estimate will be used in the proof of Lemma 2.

The following theorem is due to Kotelenez (see Theorem 2.1 and Theorem 3.4 in [20]; see also Theorem 3.1 in [19]).

THEOREM 8. *Fix $\alpha \in (0, \bar{\alpha})$. There exists a measurable function*

$$\mathcal{G}^\epsilon : \mathbb{B}_\alpha \times \mathbb{B}_0^T \to \mathcal{C}([0,T]:\mathbb{B}_\alpha)$$

*such that, for any filtered probability space $(\Omega, \mathcal{F}, \mathbb{P}, \{\mathcal{F}_t\})$ with a Brownian sheet $B$ as above and $x \in \mathbb{B}_\alpha$, $X^{\epsilon,x} \doteq \mathcal{G}^\epsilon(x, \sqrt{\varepsilon}B)$ is the unique mild solution of (5.1) (with initial condition $x$), and satisfies $\sup_{0 \leq t \leq T} \mathbb{E}\|X^{\epsilon,x}(t)\|_0^p < \infty$ for all $p \geq 0$.*

For the rest of the section we will only consider $\alpha \in (0, \bar{\alpha})$. For $f \in \mathcal{C}([0,T] : \mathbb{B}_\alpha)$, define

$$I_x(f) \doteq \inf_u \int_{[0,T] \times \mathcal{O}} u^2(s,q) \, dq \, ds, \tag{5.11}$$

where the infimum is taken over all $u \in L^2([0,T] \times \mathcal{O})$ such that

$$\begin{aligned} f(t,r) &= \int_\mathcal{O} G(t,0,r,q) x(q) \, dq \\ &\quad + \int_{[0,t] \times \mathcal{O}} G(t,s,r,q) R(s,q,f(s,q)) \, dq \, ds \\ &\quad + \int_{[0,t] \times \mathcal{O}} G(t,s,r,q) F(s,q,f(s,q)) u(s,q) \, dq \, ds. \end{aligned} \tag{5.12}$$

The following is the main result of this section.

THEOREM 9. *Let $X^{\epsilon,x}$ be as in Theorem 8. Then $I_x$ defined by (5.11) is a rate function on $\mathcal{C}([0,T] : \mathbb{B}_\alpha)$ and the family $\{I_x, x \in \mathbb{B}_\alpha\}$ of rate functions has compact level sets on compacts. Furthermore, $\{X^{\epsilon,x}\}$ satisfies the Laplace principle on $\mathcal{C}([0,T] : \mathbb{B}_\alpha)$ with the rate function $I_x$, uniformly for $x$ in compact subsets of $\mathbb{B}_\alpha$.*

REMARK 3. 1. If Assumption 4 (2) is weakened to merely the requirement that, for every $\xi \in \mathbb{B}_\alpha$, $t \mapsto \int_\mathcal{O} G(t,0,\cdot,q) \xi(q) \, dq$ is in $\mathcal{C}([0,T] : \mathbb{B}_\alpha)$, then the proof of Theorem 9 shows that, for all $x \in \mathbb{B}_\alpha$, the large deviation principle for $\{X^{\epsilon,x}\}$ on $\mathcal{C}([0,T] : \mathbb{B}_\alpha)$ holds (but not necessarily uniformly).

2. The small noise LDP for a class of reaction-diffusion SPDEs, with $\mathcal{O} = [0,1]$ and a bounded diffusion coefficient, has been studied in [25]. A difference in the conditions on the kernel $G$ in [25] is that instead of (5.6), $G$ satisfies the $L^2$ estimate in Remark 2 (3) with $\overline{\alpha} = 1/4$. One finds that the proof of Lemma 2, which is at the heart of the proof of Theorem 9, only uses the $L^2$ estimate rather than the condition (5.6). Using this observation, one can, in a straightforward manner, extend results of [25] to the case where the diffusion coefficient, instead of being bounded, satisfies the linear growth condition (5.8).

Since the proof of Theorem 9 relies on properties of the controlled process (5.3), the first step is to prove existence and uniqueness of solutions. This follows from a standard application of Girsanov's theorem.

THEOREM 10. *Let $\mathcal{G}^\epsilon$ be as in Theorem 8 and let $u \in \mathcal{P}_2^N$ for some $N \in \mathbb{N}_0$ where $\mathcal{P}_2^N$ is as defined in (4.6). For $\epsilon > 0$ and $x \in \mathbb{B}_\alpha$, define*

$$X_x^{\epsilon,u} \doteq \mathcal{G}^\epsilon(x, \sqrt{\varepsilon} B + \mathrm{Int}(u)),$$

*where $\mathrm{Int}$ is defined in (4.7). Then $X_x^{\epsilon,u}$ is the unique solution of (5.3).*

PROOF.  Fix $u \in \mathcal{P}_2^N$. Since

$$\mathbb{E}\bigg(\exp\bigg\{-\frac{1}{\sqrt{\epsilon}}\int_{[0,T]\times\mathcal{O}} u(s,q)B(dq\,ds) - \frac{1}{2\epsilon}\int_{[0,T]\times\mathcal{O}} u^2(s,q)\,dq\,ds\bigg\}\bigg) = 1,$$

the measure $\gamma^{u,\epsilon}$ defined by

$$d\gamma^{u,\epsilon} = \exp\bigg\{-\frac{1}{\sqrt{\epsilon}}\int_{[0,T]\times\mathcal{O}} u(s,q)B(dq\,ds) - \frac{1}{2\epsilon}\int_{[0,T]\times\mathcal{O}} u^2(s,q)\,dq\,ds\bigg\}\,d\mathbb{P}$$

is a probability measure on $(\Omega, \mathcal{F}, \mathbb{P})$. Furthermore, $\gamma^{u,\epsilon}$ is mutually absolutely continuous with respect to $\mathbb{P}$ and by Girsanov's theorem (see [7], Theorem 10.14), the process $\widetilde{B} = B + \epsilon^{-1/2}\operatorname{Int}(u)$ on $(\Omega, \mathcal{F}, \gamma^{u,\epsilon}, \{\mathcal{F}_t\})$ is a Brownian sheet. Thus, by Theorem 8, $X_x^{\epsilon,u} = \mathcal{G}^\epsilon(x, \sqrt{\varepsilon}B + \operatorname{Int}(u))$ is the unique solution of (5.2), with $B$ there replaced by $\widetilde{B}$, on $(\Omega, \mathcal{F}, \gamma^{u,\epsilon}, \{\mathcal{F}_t\})$. However, equation (5.2) with $\widetilde{B}$ is precisely the same as equation (5.3), and since $\gamma^{u,\epsilon}$ and $\mathbb{P}$ are mutually absolutely continuous, we get that $X_x^{\epsilon,u}$ is the unique solution of (5.3) on $(\Omega, \mathcal{F}, \mathbb{P}, \{\mathcal{F}_t\})$. This completes the proof. □

In the next subsection we will study, under the standing assumption of this section, the following two basic qualitative results regarding the processes $X_x^{\epsilon,u}$. The first is simply the controlled, zero-noise version of the theorem just stated and its proof, being very similar to the proof of Theorem 8, is omitted. The next is a standard convergence result whose proof is given in Section 5.2.

THEOREM 11.  *Fix $x \in \mathbb{B}_\alpha$ and $u \in L^2([0,T]\times\mathcal{O})$. Then there is a unique function $f$ in $\mathcal{C}([0,T]:\mathbb{B}_\alpha)$ which satisfies equation (5.12).*

In analogy with the notation $X_x^{\varepsilon,u}$ for the solution of (5.3), we will denote the unique solution $f$ given by Theorem 11 as $X_x^{0,u}$. Let $\theta:[0,1)\to[0,1)$ be a measurable map such that $\theta(r) \to \theta(0) = 0$ as $r \to 0$.

THEOREM 12.  *Let $M < \infty$, and suppose that $x^\epsilon \to x$ and $u^\epsilon \to u$ in distribution as $\epsilon \to 0$ with $\{u^\epsilon\} \subset \mathcal{P}_2^M$. Then $X_{x^\epsilon}^{\theta(\epsilon),u^\epsilon} \to X_x^{0,u}$ in distribution.*

PROOF OF THEOREM 9.  Define the map $\mathcal{G}^0:\mathbb{B}_\alpha \times \mathbb{B}_0^T \to \mathcal{C}([0,T]:\mathbb{B}_\alpha)$ as follows. For $x \in \mathbb{B}_\alpha$ and $\phi \in \mathbb{B}_0^T$ of the form $\phi(t,x) \doteq \operatorname{Int}(u)(t,x)$ for some $u \in L^2([0,T]\times\mathcal{O})$, we define $\mathcal{G}^0(x,\phi) = X_x^{0,u}$. Set $\mathcal{G}^0(x,\phi) = 0$ for all other $\phi \in \mathbb{B}_0^T$. In view of Theorem 7, it suffices to show that $(\mathcal{G}^\epsilon, \mathcal{G}^0)$ satisfy Assumption 3 with $\mathcal{E}_0$ and $\mathcal{E}$ there replaced by $\mathbb{B}_\alpha$ and $\mathcal{C}([0,T];\mathbb{B}_\alpha)$ respectively; and for all $f \in \mathcal{E}$, the map $x \mapsto I_x(f)$ is l.s.c. The latter property and the first part of Assumption 3 is immediate on applying Theorems 11 and 12 with $\theta = 0$. The second part of Assumption 3 follows on applying Theorem 12 with $\theta(r) = r$, $r \in [0,1)$. □

5.2. *Qualitative properties of controlled stochastic reaction-diffusion equations.* This section is devoted to the proof of Theorem 12. Our first result shows that $L^p$ bounds hold for controlled SDEs, uniformly when the initial condition and controls lie in compact sets and $\varepsilon \in [0,1)$. Note, in particular, that $\varepsilon = 0$ is allowed.

LEMMA 1. *If $K$ is any compact subset of $\mathbb{B}_\alpha$ and $M < \infty$, then for all $p \in [1, \infty)$,*

$$\sup_{u \in \mathcal{P}_2^M} \sup_{x \in K} \sup_{\epsilon \in [0,1)} \sup_{(t,r) \in [0,T] \times \mathcal{O}} \mathbb{E}|X_x^{\epsilon,u}(t,r)|^p < \infty.$$

PROOF. By Doob's inequality, there exists a suitable constant $c_1$ such that

$$\mathbb{E}|X_x^{\epsilon,u}(t,r)|^p$$
$$\leq c_1 \left| \int_\mathcal{O} G(t,0,r,q) x(q) \, dq \right|^p$$
$$+ c_1 \mathbb{E} \left| \int_0^t \int_\mathcal{O} G(t,s,r,q) R(s,q,X_x^{\epsilon,u}(s,q)) \, dq \, ds \right|^p$$
$$+ c_1 \mathbb{E} \left[ \int_0^t \int_\mathcal{O} |G(t,s,r,q)|^2 |F(s,q,X_x^{\epsilon,u}(s,q))|^2 \, dq \, ds \right]^{p/2}$$
$$+ c_1 \mathbb{E} \left[ \int_0^t \int_\mathcal{O} |G(t,s,r,q)| |F(s,q,X_x^{\epsilon,u}(s,q))| |u(s,q)| \, dq \, ds \right]^p.$$

Using (5.8) and the Cauchy–Schwarz inequality, the right-hand side above can be bounded by

$$c_2 \left[ 1 + \mathbb{E} \left( \int_0^t \int_\mathcal{O} |G(t,s,r,q)|^2 |X_x^{\epsilon,u}(s,q)|^2 \, dq \, ds \right)^{p/2} \right].$$

Hölder's inequality yields for $p > 2$ that

$$\Lambda_p(t) \leq c_2 \left[ 1 + \left( \int_0^T \int_\mathcal{O} |G(t,s,r,q)|^{2\tilde{p}} \, dq \, ds \right)^{(p-2)/2} \int_0^t \Lambda_p(s) \, ds \right],$$

where $\Lambda_p(t) = \sup_{u \in P_2^M} \sup_{x \in K} \sup_{\epsilon \in [0,1)} \sup_{r \in \mathcal{O}} \mathbb{E}|X_x^{\epsilon,u}(t,r)|^p$ and $\tilde{p} = \frac{p}{p-2}$. Choose $p_0$ large enough that $(\frac{2p_0}{p_0-2} - 1)(1 - 2\bar{\alpha}) < 1$. Using (5.4) and (5.5), we have for all $p \geq p_0$ that

$$\left[ \int_0^T \int_\mathcal{O} |G(t,s,r,q)|^{2\tilde{p}} \, dq \, ds \right]^{(p-2)/2} \leq c_3 T^{(1-(2\tilde{p}-1)(1-2\bar{\alpha}))(p-2)/2}.$$

Thus, for every $p \geq p_0$, there exists a constant $c_4$ such that $\Lambda_p(t) \leq c_4[1 + \int_0^t \Lambda_p(s) \, ds]$. The result now follows from Gronwall's lemma. □

The following lemma will be instrumental in proving tightness and weak convergence in Banach spaces such as $\mathbb{B}_\alpha$ and $\mathbb{B}_\alpha^T$.

LEMMA 2. *Let $\mathcal{A} \subseteq \mathcal{P}_2$ be a family such that, for all $p \geq 2$,*

$$\sup_{f \in \mathcal{A}} \sup_{(t,r) \in [0,T] \times \mathcal{O}} \mathbb{E}|f(t,r)|^p < \infty. \tag{5.13}$$

*Also, let $\mathcal{B} \subseteq \mathcal{P}_2^M$ for some $M < \infty$. For $f \in \mathcal{A}$ and $u \in \mathcal{B}$, define*

$$\Psi_1(t,r) \doteq \int_0^t \int_\mathcal{O} G(t,s,r,q) f(s,q) B(dq\,ds),$$

$$\Psi_2(t,r) \doteq \int_0^t \int_\mathcal{O} G(t,s,r,q) f(s,q) u(s,q)\, dq\,ds,$$

*where the dependence on $f$ and $u$ is not made explicit in the notation. Then for any $\alpha < \bar{\alpha}$ and $i = 1, 2$,*

$$\sup_{f \in \mathcal{A}, u \in \mathcal{B}} \mathbb{E}\left\{ \sup_{\rho((t,r),(s,q))<1} \frac{|\Psi_i(t,r) - \Psi_i(s,q)|}{\rho((t,r),(s,q))^\alpha} \right\} < \infty.$$

PROOF. We will prove the result for $i = 1$; the proof for $i = 2$ is identical (except an additional application of the Cauchy–Schwarz inequality) and, thus, it is omitted. Henceforth, we write, for simplicity, $\Psi_1$ as $\Psi$. We will apply Theorem 6 of [15], according to which it suffices to show that, for all $0 \leq t_1 < t_2 \leq T$, $r_1, r_2 \in \mathcal{O}$,

$$\sup_{f \in \mathcal{A}, u \in \mathcal{B}} \mathbb{E}|\Psi(t_2, r_2) - \Psi(t_1, r_1)|^p \leq c_p (\hat{\omega}(\rho((t_1, r_1), (t_2, r_2))))^p, \tag{5.14}$$

for a suitable constant $c_p$; a $p > 2$; and a function $\hat{\omega}: [0, \infty) \to [0, \infty)$ satisfying

$$\int_0^1 \frac{\hat{\omega}(u)}{u^{1+\alpha+(d+1)/p}}\, du < \infty.$$

We will show that (5.14) holds with $\hat{\omega}(u) = u^{\alpha_0}$ for some $\alpha_0 \in (\alpha, \bar{\alpha})$ and all $p$ sufficiently large. This will establish the result.

Fix $\alpha_0, \tilde{\alpha}$ such that $\alpha < \alpha_0 < \tilde{\alpha} < \bar{\alpha}$ and let $t_1 < t_2$, $r_1, r_2 \in \mathcal{O}$ and $p > 2$. We will need $p$ to be sufficiently large and the choice of $p$ will be fixed in the course of the proof. By Doob's inequality, there exists a constant $c_1$ such that

$$\begin{aligned}
\mathbb{E}&|\Psi(t_2, r_2) - \Psi(t_1, r_1)|^p \\
&\leq c_1 \mathbb{E}\left[ \int_0^T \int_\mathcal{O} |G(t_2, s, r_2, q) - G(t_1, s, r_1, q)|^2 |f(s,q)|^2\, dq\,ds \right]^{p/2}.
\end{aligned} \tag{5.15}$$

Let $\tilde{p} = p/(p-2)$ and $\delta = 4/p$. Note that $(2-\delta)\tilde{p} = \delta p/2 = 2$. Hölder's inequality, (5.9) and (5.13) give that the right-hand side of (5.15) is bounded by

$$
\begin{aligned}
(5.16) \quad & c_1 \left[ \int_0^T \int_{\mathcal{O}} |G(t_2, s, r_2, q) - G(t_1, s, r_1, q)|^{(2-\delta)\tilde{p}} \, dq \, ds \right]^{(p-2)/2} \\
& \times \left[ \int_0^T \int_{\mathcal{O}} |G(t_2, s, r_2, q) - G(t_1, s, r_1, q)|^{\delta p/2} \mathbb{E}|f(s,q)|^p \, dq \, ds \right] \\
& \leq c_2 \left[ \int_0^T \int_{\mathcal{O}} |G(t_2, s, r_2, q) - G(t_1, s, r_1, q)|^2 \, dq \, ds \right]^{(p-2)/2}
\end{aligned}
$$

for a suitable constant $c_2$ that is independent of $f$. From Remark 2(3), the expression in (5.16) can be bounded (for $p$ large enough) by

$$c_3 \rho((t_1, r_1), (t_2, r_2))^{\tilde{\alpha}(p-2)} \leq c_4 \rho((t_1, r_1), (t_2, r_2))^{\alpha_0 p}.$$

The result follows. □

The next result will be used to prove the stochastic integral in (5.3) converges to 0 in $\mathcal{C}([0,T] \times \mathcal{O})$, which will be strengthened shortly.

LEMMA 3. *Let $\mathcal{A}$ and $\Psi_1$ be as in Lemma 2 and let $Z_f^\epsilon \doteq \sqrt{\epsilon} \Psi_1$. Then for every sequence $\{f^\varepsilon\} \subset \mathcal{A}$, $Z_{f^\varepsilon}^\epsilon \xrightarrow{\mathbb{P}} 0$ in $\mathcal{C}([0,T] \times \mathcal{O})$, as $\varepsilon \to 0$.*

PROOF. Arguments similar to those lead to (5.16) along with (5.4), (5.5) and (5.13) yielding that $\sup_{f \in \mathcal{A}} \mathbb{E}|\Psi_1(t,r)|^2 < \infty$. This shows that for each $(t,r) \in [0,T] \times \mathcal{O}$, $Z_{f^\varepsilon}^\epsilon(t,r) \xrightarrow{\mathbb{P}} 0$ (in fact in $L^2$). Defining

$$\omega(x, \delta) \doteq \sup\{|x(t,r) - x(t',r')| : \rho((t,r),(t',r')) \leq \delta\}$$

for $x \in \mathcal{C}([0,T] \times \mathcal{O})$ and $\delta \in (0,1)$, we see that $\omega(Z_{f^\varepsilon}^\epsilon, \delta) = \sqrt{\epsilon} \delta^\alpha M_{f^\varepsilon}^\epsilon$, where $M_f^\epsilon \doteq \sup_{\rho((t,r),(s,q)) < 1} \frac{|\Psi_1(t,r) - \Psi_1(s,q)|}{\rho((t,r),(s,q))^\alpha}$. Therefore, from Lemma 2,

$$\lim_{\delta \to 0} \lim_{\epsilon \to 0} \mathbb{E} \omega(Z_{f^\varepsilon}^\epsilon, \delta) = 0.$$

The result now follows from Theorem 14.5 of [17]. □

We now establish the main convergence result.

PROOF OF THEOREM 12. Given $x \in K, u \in P_2^M, \epsilon \in [0,1)$, define

$$Z_{1,x}^{\epsilon,u}(t,r) = \int_{\mathcal{O}} G(t,0,r,q) x(q) \, dq,$$

$$Z_{2,x}^{\epsilon,u}(t,r) = \int_0^t \int_{\mathcal{O}} G(t,s,r,q) R(s,q,X_x^{\theta(\epsilon),u}(s,q))\, dq\, ds,$$

$$Z_{3,x}^{\epsilon,u}(t,r) = \sqrt{\theta(\epsilon)} \int_0^t \int_{\mathcal{O}} G(t,s,r,q) F(s,q,X_x^{\theta(\epsilon),u}(s,q)) B(dq\, ds),$$

$$Z_{4,x}^{\epsilon,u}(t,r) = \int_0^t \int_{\mathcal{O}} G(t,s,r,q) F(s,q,X_x^{\theta(\epsilon),u}(s,q)) u(s,q)\, dq\, ds.$$

We first show that each $Z_{i,x^\epsilon}^{\epsilon,u^\epsilon}$ is tight in $\mathcal{C}([0,T]:\mathbb{B}_\alpha)$, for $i=1,2,3,4$. For $i=1$, this follows from part 2 of Assumption 4. Recalling that $\mathbb{B}_{\alpha^*}^T$ is compactly embedded in $\mathbb{B}_\alpha^T$ for $\bar{\alpha} > \alpha^* > \alpha$, it suffices to show that, for some $\alpha^* \in (\alpha, \bar{\alpha})$,

(5.17) $$\sup_{\epsilon \in (0,1)} \mathbb{P}[\|Z_{i,x^\epsilon}^{\epsilon,u^\epsilon}\|_{\mathbb{B}_{\alpha^*}^T} > K] \to 0 \quad \text{as } K \to \infty \text{ for } i=2,3,4.$$

For $i=2,4$, (5.17) is an immediate consequence of

$$\sup_{\epsilon \in (0,1)} \mathbb{E}\|Z_{i,x^\epsilon}^{\epsilon,u^\epsilon}\|_{\mathbb{B}_{\alpha^*}^T} < \infty,$$

as follows from Lemma 2, the linear growth condition (5.8) and Lemma 1. For $i=3$, in view of Lemma 3, it suffices to establish

$$\sup_{\epsilon \in (0,1)} \mathbb{E}[Z_{3,x^\epsilon}^{\epsilon,u^\epsilon}]_{\mathbb{B}_{\alpha^*}^T} < \infty,$$

where for $z \in \mathbb{B}_\alpha^T$, $[z]_{\mathbb{B}_\alpha^T} = \|z\|_{\mathbb{B}_\alpha^T} - \|z\|_0$. Once more, this follows as an immediate consequence of Lemma 2, the linear growth condition (5.8) and Lemma 1.

Having shown tightness of $Z_{i,x^\epsilon}^{\epsilon,u^\epsilon}$ for $i=1,2,3,4$, we can extract a subsequence along which each of these processes and $X_{x^\epsilon}^{\epsilon,u^\epsilon}$ converges in distribution in $\mathcal{C}([0,T]:\mathbb{B}_\alpha)$. Let $Z_{i,x}^{0,u}$ and $X_x^{0,u}$ denote the respective limits. We will show that

$$Z_{1,x}^{0,u}(t,r) = \int_{\mathcal{O}} G(t,0,r,q) x(q)\, dq,$$

$$Z_{2,x}^{0,u}(t,r) = \int_0^t \int_{\mathcal{O}} G(t,s,r,q) R(s,q,X_x^{0,u}(s,q))\, dq\, ds,$$

(5.18)

$$Z_{3,x}^{0,u}(t,r) = 0,$$

$$Z_{4,x}^{0,u}(t,r) = \int_0^t \int_{\mathcal{O}} G(t,s,r,q) F(s,q,X_x^{0,u}(s,q)) u(s,q)\, dq\, ds.$$

The uniqueness result Theorem 11 will then complete the proof.

Convergence for $i=1$ follows from part 2 of Assumption 4. The case $i=3$ follows from Lemma 3, Lemma 1 and the linear growth condition. To deal with the cases $i=2,4$, we invoke the Skorokhod Representation Theorem [21], which allows us to assume with probability one convergence for the purposes of identifying the limits. We give the proof of convergence only for the harder case $i = 4$. Denote the right-hand side of (5.18) by $\hat{Z}^{0,u}_{4,x}(t,r)$. Then

$$
\begin{aligned}
|Z^{\epsilon,u^\epsilon}_{4,x^\epsilon}&(t,r) - \hat{Z}^{0,u}_{4,x}(t,r)| \\
&\leq \int_0^t \int_{\mathcal{O}} |G(t,s,r,q)| |F(s,q,X^{\epsilon,u^\epsilon}_{x^\epsilon}(s,q)) \\
&\qquad\qquad - F(s,q,X^{0,u}_x(s,q))| |u^\epsilon(s,q)| \, dq \, ds \\
&\quad + \left| \int_0^t \int_{\mathcal{O}} G(t,s,r,q) F(s,q,X^{0,u}_x(s,q))(u^\epsilon(s,q) - u(s,q)) \, dq \, ds \right|.
\end{aligned}
\tag{5.19}
$$

By the Cauchy–Schwarz inequality, equation (5.9) and the uniform Lipschitz property of $F$ we see that, for a suitable constant $c \in (0, \infty)$, the first term on the right-hand side of (5.19) can be bounded above by

$$
\sqrt{M} \left[ \int_0^t \int_{\mathcal{O}} |G(t,s,r,q)|^2 |F(s,q,X^{\epsilon,u^\epsilon}_{x^\epsilon}(s,q)) - F(s,q,X^{0,u}_x(s,q))|^2 \, dq \, ds \right]^{1/2}
$$

$$
\leq c \left( \sup_{(s,q) \in [0,T] \times \mathcal{O}} |X^{\epsilon,u^\epsilon}_{x^\epsilon}(s,q) - X^{0,u}_x(s,q)| \right),
$$

and thus converges to 0 as $\varepsilon \to 0$. The second term in (5.19) converges to 0 as well, since $u^\varepsilon \to u$ and

$$
\int_0^t \int_{\mathcal{O}} (G(t,s,r,q) F(s,q,X^{0,u}_x(s,q)))^2 \, dq \, ds < \infty.
$$

By uniqueness of limits and noting that $\hat{Z}^{0,u}_{4,x}$ is a continuous random field, we see that $Z^{0,u}_{4,x} = \hat{Z}^{0,u}_{4,x}$ and the proof is complete. $\square$

**6. Other infinite dimensional models.** The key ingredients in the proof of the LDP for the solution of the infinite dimensional SDE studied in Section 5 are the qualitative properties in Theorems 11 and 12 of the controlled SDE (5.3). Once these properties are verified, the LDP follows as an immediate consequence of Theorem 7. Furthermore, one finds that the estimates needed for the proof of Theorems 11 and 12 are essentially the same as those needed for establishing unique solvability of (5.1). This is a common theme that appears in all proofs of LDPs, for small noise stochastic dynamical systems, that are based on variational representations such as in Section 3. Indeed, one can argue that the variational representation approach makes the small noise large deviation analysis a transparent and a largely straightforward exercise, once one has the estimates for the unique solvability of the stochastic equation. This statement has been affirmed by several recent works on Freidlin–Wentzell large deviations for infinite dimensional SDEs that are based on the variational representation approach (specifically Theorem 2), and carry out the verification of statements analogous to Theorems 11 and 12. Some of these works are summarized below.

6.1. *SDEs driven by infinitely many Brownian motions.* Ren and Zhang [23] consider a SDE driven by infinitely many Brownian motions with non-Lipschitz diffusion coefficients. Prior results on strong existence and uniqueness of the solutions to the SDE yield continuous (in time and initial condition) random field solutions. The authors prove a small noise LDP in the space $C([0,T] \times \mathbb{R}^d)$. The proof relies on the representation formula for an infinite sequence of real Brownian motions $\{\beta_i\}$ given in Theorem 2 and the general Laplace principle of the form in Theorem 6. Non-Lipschitz coefficients make the standard discretization and approximation approach intractable for this example. The authors verify the analogues of Theorems 11 and 12 in Theorem 3.1, Lemmas 3.4 and 3.11 of the cited paper. In the final section of the paper, Schilder's theorem for Brownian motion on the group of homeomorphisms of the circle is obtained. The proof here is also by verifying of steps analogous to Theorems 11 and 12 regarding solvability and convergence in the space of homeomorphisms. Once more, exponential probability estimates with the natural metric on the space of homeomorphisms, needed in the standard proofs of the LDP, do not appear to be straightforward. Using similar ideas based on representations for infinite dimensional Brownian motions, a LDP for flows of homeomorphisms, extending results of the final section of [23] to multi-dimensional SDES with nonLipschitz coefficients, has been studied in [24].

6.2. *Stochastic PDE with varying boundary conditions.* Wang and Duan [27] study stochastic parabolic PDEs with rapidly varying random dynamical boundary conditions. The formulation of the SPDE as an abstract stochastic evolution equation in an appropriate Hilbert space leads to a non-Lipschitz nonlinearity with polynomial growth. Deviations of the solution from the limiting effective system (as the parameter governing the rapid component approaches its limit) are studied by establishing a large deviation principle. The proof of the LDP uses the variational representation for functionals of a Hilbert space valued Wiener process as in Theorem 3 and the general Laplace principle given in Theorem 5. Once more, the hardest part in the analysis is establishing the well-posedness (i.e., existence, uniqueness) of the stochastic evolution equation. Once estimates for existence/uniqueness are available, the proof of the LDP becomes a straightforward verification of Assumption 1.

6.3. *Stochastic Navier–Stokes equation.* Sritharan and Sundar [26] study small noise large deviations for a two-dimensional Navier–Stokes equation in an (possibly) unbounded domain and with multiplicative noise. The equation can be written as an abstract stochastic evolution equation in an appropriate function space. The solution lies in the Polish space $C([0,T]:H) \cap L^2([0,T]:V)$ for some Hilbert spaces $H$ and $V$ and can be expressed as $\mathcal{G}^\varepsilon(\sqrt{\varepsilon}W)$ for a $H$ valued Wiener process $W$. Authors prove existence and uniqueness of solutions and then apply Theorem 5 by verifying Assumption 1 for their model.

## APPENDIX

PROOF OF THEOREM 5. For the first part of the theorem, we need to show that, for all compact subsets $K$ of $\mathcal{E}_0$ and each $M < \infty$,

$$\Lambda_{M,K} \doteq \bigcup_{x \in K} \{f \in \mathcal{E} : I_x(f) \leq M\}$$

is a compact subset of $\mathcal{E}$. To establish this, we will show that $\Lambda_{M,K} = \bigcap_{n \geq 1} \Gamma_{2M+1/n,K}$. In view of Assumption 1, the compactness of $\Lambda_{M,K}$ will then follow. Let $f \in \Lambda_{M,K}$. There exists $x \in K$ such that $I_x(f) \leq M$. We can now find, for each $n \geq 1$, $u_n \in L^2([0,T]:H_0)$ such that $f = \mathcal{G}^0(x, \int_0^\cdot u_n(s) \, ds)$ and $\frac{1}{2} \int_0^T \|u_n(s)\|_0^2 \, ds \leq M + \frac{1}{2n}$. In particular, $u_n \in S^{2M+1/n}(H_0)$, and so $f \in \Gamma_{2M+1/n,K}$. Since $n \geq 1$ is arbitrary, we have $\Lambda_{M,K} \subseteq \bigcap_{n \geq 1} \Gamma_{2M+1/n,K}$. Conversely, suppose $f \in \Gamma_{2M+1/n,K}$, for all $n \geq 1$. Then, for every $n \geq 1$, there exists $x_n \in K, u_n \in S^{2M+1/n}$ such that $f = \mathcal{G}^0(x_n, \int_0^\cdot u_n(s) \, ds)$. In particular, we have $I_{x_n}(f) \leq M + \frac{1}{2n}$. Recalling that the map $x \mapsto I_x(f)$ is l.s.c. and $K$ is compact, we now see on sending $n \to \infty$ that, for some $x \in K$, $I_x(f) \leq M$. Thus, $f \in \Lambda_{M,K}$ and the inclusion $\bigcap_{n \geq 1} \Gamma_{2M+1/n,K} \subseteq \Lambda_{M,K}$ follows. This proves the first part of the theorem.

For the second part of the theorem, consider an $x \in \mathcal{E}_0$ and let $\{x^\epsilon, \varepsilon > 0\} \subseteq \mathcal{E}_0$ be such that $x^\epsilon \to x$ as $\epsilon \to 0$. Fix a bounded and continuous function $h : \mathcal{E} \to \mathbb{R}$. It suffices to show (2.1) (upper bound) and (2.2) (lower bound), with $X^\varepsilon$ there replaced by $X^{\varepsilon,x^\varepsilon}$ and $I$ replaced by $I_x$. For notational convenience, we will write $\mathcal{P}_2(H_0), \mathcal{P}_2^N(H_0), S^N(H_0)$ simply as $\mathcal{P}_2, \mathcal{P}_2^N, S^N$ respectively.

**Proof of the upper bound.** From Theorem 3,

$$\begin{aligned}(7.1) \quad &-\epsilon \log \mathbb{E}\left[\exp\left(-\frac{1}{\epsilon} h(X^{\epsilon,x^\varepsilon})\right)\right] \\ &= \inf_{u \in \mathcal{P}_2} \mathbb{E}\left[\frac{1}{2} \int_0^T \|u(s)\|_0^2 \, ds + h \circ \mathcal{G}^\epsilon\left(x^\varepsilon, \sqrt{\varepsilon}W + \int_0^\cdot u(s) \, ds\right)\right].\end{aligned}$$

Fix $\delta \in (0,1)$. Then for every $\epsilon > 0$ there exists $u^\epsilon \in \mathcal{P}_2$ such that the right-hand side of (7.1) is bounded below by

$$\mathbb{E}\left[\tfrac{1}{2}\int_0^T \|u^\epsilon(s)\|_0^2\,ds + h \circ \mathcal{G}^\epsilon\left(x^\epsilon, \sqrt{\epsilon}W + \int_0^\cdot u^\epsilon(s)\,ds\right)\right] - \delta.$$

Using the fact that $h$ is bounded, we can assume without loss of generality (we refer the reader to the proof of Theorem 4.4 of [3] where a similar argument is used) that, for some $N \in (0,\infty)$,

$$\sup_{\epsilon > 0} \int_0^T \|u^\epsilon(s)\|_0^2\,ds \leq N \qquad \text{a.s.}$$

In order to prove the upper bound, it suffices to show that

$$\liminf_{\epsilon \to 0} \mathbb{E}\left[\tfrac{1}{2}\int_0^T \|u^\epsilon(s)\|_0^2\,ds + h \circ \mathcal{G}^\epsilon\left(x^\epsilon, \sqrt{\epsilon}W + \int_0^\cdot u^\epsilon(s)\,ds\right)\right]$$
$$\geq \inf_{f \in \mathcal{E}}\{I_x(f) + h(f)\}.$$

Pick a subsequence (relabeled by $\epsilon$) along which $u^\epsilon$ converges in distribution to some $u \in \mathcal{P}_2^N$ (as $S^N$ valued random variables). We now infer from the Assumption 1 that

$$\liminf_{\epsilon \to 0} \mathbb{E}\left[\tfrac{1}{2}\int_0^T \|u^\epsilon(s)\|_0^2\,ds + h \circ \mathcal{G}^\epsilon\left(x^\epsilon, \sqrt{\epsilon}W + \int_0^\cdot u^\epsilon(s)\,ds\right)\right]$$
$$\geq \mathbb{E}\left[\tfrac{1}{2}\int_0^T \|u(s)\|_0^2\,ds + h \circ \mathcal{G}^0\left(x, \int_0^\cdot u(s)\,ds\right)\right]$$
$$\geq \inf_{\{(f,u) \in \mathcal{E} \times L^2([0,T]:H_0) : f = \mathcal{G}^0(x, \int_0^\cdot u(s)\,ds)\}} \left\{\tfrac{1}{2}\int_0^T \|u(s)\|_0^2\,ds + h(f)\right\}$$
$$\geq \inf_{f \in \mathcal{E}}\{I_x(f) + h(f)\}.$$

**Proof of the lower bound.** We need to show that

$$\limsup_{\epsilon \to 0} -\epsilon \log \mathbb{E}\left[\exp\left(-\tfrac{1}{\epsilon}h(X^{\epsilon,x^\epsilon})\right)\right] \leq \inf_{f \in \mathcal{E}}\{I_x(f) + h(f)\}.$$

Without loss of generality, we can assume that $\inf_{f \in \mathcal{E}}\{I_x(f) + h(f)\} < \infty$. Let $\delta > 0$ be arbitrary, and let $f_0 \in \mathcal{E}$ be such that

(7.2) $$I_x(f_0) + h(f_0) \leq \inf_{f \in \mathcal{E}}\{I_x(f) + h(f)\} + \frac{\delta}{2}.$$

Choose $\tilde{u} \in L^2([0,T]:H_0)$ such that,

(7.3) $$\tfrac{1}{2}\int_0^T \|\tilde{u}(s)\|_0^2\,ds \leq I_x(f_0) + \frac{\delta}{2} \quad \text{and} \quad f_0 = \mathcal{G}^0\left(x, \int_0^\cdot \tilde{u}(s)\,ds\right).$$

Then, from Theorem 3,

$$\limsup_{\epsilon \to 0} -\epsilon \log \mathbb{E}\left[\exp\left(-\frac{1}{\epsilon} h(X^{\epsilon,x^\varepsilon})\right)\right]$$

(7.4)
$$= \limsup_{\epsilon \to 0} \inf_{u \in \mathcal{A}} \mathbb{E}\left[\frac{1}{2}\int_0^T \|u(s)\|_0^2\, ds + h \circ \mathcal{G}^\epsilon\left(x^\varepsilon, \sqrt{\epsilon}W + \int_0^\cdot u(s)\, ds\right)\right]$$

$$\le \limsup_{\epsilon \to 0} \mathbb{E}\left[\frac{1}{2}\int_0^T \|\tilde{u}(s)\|_0^2\, ds + h \circ \mathcal{G}^\epsilon\left(x^\varepsilon, \sqrt{\epsilon}W + \int_0^\cdot \tilde{u}(s)\, ds\right)\right]$$

$$= \frac{1}{2}\int_0^T \|\tilde{u}(s)\|_0^2\, ds + \limsup_{\epsilon \to 0} \mathbb{E}\left[h \circ \mathcal{G}^\epsilon\left(x^\varepsilon, \sqrt{\epsilon}W + \int_0^\cdot \tilde{u}(s)\, ds\right)\right].$$

By Assumption 1, $\lim_{\epsilon \to 0} \mathbb{E}[h \circ \mathcal{G}^\epsilon(x^\varepsilon, \sqrt{\epsilon}W + \int_0^\cdot \tilde{u}(s)\, ds)] = h(\mathcal{G}^0(x, \int_0^\cdot \tilde{u}(s)\, ds)) = h(f_0)$. Thus, in view of (7.2) and (7.3), the expression (7.4) can be at most $\inf_{f \in \mathcal{E}}\{I(f) + h(f)\} + \delta$. Since $\delta$ is arbitrary, the proof is complete. $\square$

PROOF OF THEOREM 6. From Remark 1, we can regard $\beta$ as an $H$ valued $Q$-Wiener process, where $H = \bar{l}_2$ and $Q$ is a trace class operator, as defined in Remark 1. Also, one can check that $H_0 \doteq Q^{1/2}H = l_2$. Since the embedding map $i: C([0,T]: \bar{l}_2) \to C([0,T]: \mathbb{R}^\infty)$ is measurable (in fact, continuous), $\hat{\mathcal{G}}^\varepsilon: \mathcal{E}_0 \times C([0,T]: \bar{l}_2) \to \mathcal{E}$ defined as $\hat{\mathcal{G}}^\varepsilon(x, \sqrt{\varepsilon}v) \doteq \mathcal{G}^\varepsilon(x, \sqrt{\varepsilon}i(v))$, $(x,v) \in \mathcal{E}_0 \times C([0,T]: \bar{l}_2)$ is a measurable map for every $\varepsilon \ge 0$. Note also that, for $\varepsilon > 0$, $X^{\varepsilon,x} = \hat{\mathcal{G}}^\varepsilon(x, \sqrt{\varepsilon}\beta)$ a.s. Since Assumption 2 holds, we have that 1 and 2 of Assumption 1 are satisfied with $\mathcal{G}^\varepsilon$ there replaced by $\hat{\mathcal{G}}^\varepsilon$ for $\varepsilon \ge 0$ and $W$ replaced with $\beta$. Define $\hat{I}_x(f)$ by the right-hand side of (4.3) with $\mathcal{G}^0$ replaced by $\hat{\mathcal{G}}^0$. Clearly, $I_x(f) = \hat{I}_x(f)$ for all $(x,f) \in \mathcal{E}_0 \times \mathcal{E}$. The result is now an immediate consequence of Theorem 5. $\square$

PROOF OF THEOREM 7. Let $\{\phi_i\}_{i=1}^\infty$ be a CONS in $L^2(\mathcal{O})$ and let

$$\beta_i(t) \doteq \int_{[0,t] \times \mathcal{O}} \phi_i(x) B(ds\, dx), \qquad t \in [0,T], i = 1, 2, \dots.$$

Then $\beta \equiv \{\beta_i\}$ is a sequence of independent standard real Brownian motions and can be regarded as a $(S, \mathcal{S})$ valued random variable. Furthermore, (2.5) is satisfied and from Proposition 3, there is a measurable map $g: C([0,T]: \mathbb{R}^\infty) \to C([0,T] \times \mathcal{O}: \mathbb{R})$ such that $g(\beta) = B$ a.s. Define, for $\varepsilon > 0$, $\hat{\mathcal{G}}^\varepsilon: \mathcal{E}_0 \times C([0,T]: \mathbb{R}^\infty) \to \mathcal{E}$ as $\hat{\mathcal{G}}^\varepsilon(x, \sqrt{\varepsilon}v) \doteq \mathcal{G}^\varepsilon(x, \sqrt{\varepsilon}g(v))$, $(x,v) \in \mathcal{E}_0 \times C([0,T]: \mathbb{R}^\infty)$. Clearly, $\hat{\mathcal{G}}^\varepsilon$ is a measurable map and $\hat{\mathcal{G}}^\varepsilon(x, \sqrt{\varepsilon}\beta) = X^{\varepsilon,x}$ a.s. Next, note that

$$S_{ac} \doteq \left\{v \in C([0,T]: \mathbb{R}^\infty): v(t) = \int_0^t \hat{u}(s)\, ds,\right.$$

$$t \in [0,T], \text{ for some } \hat{u} \in L^2([0,T]:l^2)\Big\}$$

is a measurable subset of $S$. For $\hat{u} \in L^2([0,T]:l_2)$, define $u_{\hat{u}} \in L^2([0,T] \times \mathcal{O})$ as

$$u_{\hat{u}}(t,x) = \sum_{i=1}^{\infty} \hat{u}_i(t)\phi_i(x), \qquad (t,x) \in [0,T] \times \mathcal{O}.$$

Define $\hat{\mathcal{G}}^0 : \mathcal{E}_0 \times C([0,T]:\mathbb{R}^\infty) \to \mathcal{E}$ as

$$\hat{\mathcal{G}}^0(x,v) \doteq \mathcal{G}^0(x, \text{Int}(u_{\hat{u}})) \qquad \text{if } v = \int_0^\cdot \hat{u}(s)\,ds \text{ and } \hat{u} \in L^2([0,T]:l_2).$$

We set $\hat{\mathcal{G}}^0(x,v) = 0$ for all other $(x,v)$. Note that

$$\left\{\hat{\mathcal{G}}^0\left(x, \int_0^\cdot \hat{u}(s)\,ds\right) : \hat{u} \in S^M(l_2), x \in K\right\}$$
$$= \{\mathcal{G}^0(x, \text{Int}(u)) : u \in S^M, x \in K\}.$$

Since Assumption 3 holds, we have that 1 of Assumption 2 holds with $\mathcal{G}^0$ there replaced by $\hat{\mathcal{G}}^0$. Next, an application of Girsanov's theorem gives that, for every $\hat{u}^\varepsilon \in \mathcal{P}_2^M(l_2)$,

$$g\left(\beta + \frac{1}{\sqrt{\epsilon}} \int_0^\cdot \hat{u}^\epsilon(s)\,ds\right) = B + \frac{1}{\sqrt{\epsilon}} \text{Int}(u_{\hat{u}^\varepsilon}) \qquad \text{a.s.}$$

In particular, for every $M < \infty$ and families $\{\hat{u}^\epsilon\} \subset \mathcal{P}_2^M(l_2)$ and $\{x^\varepsilon\} \subset \mathcal{E}_0$, such that $\hat{u}^\epsilon$ converges in distribution [as $S^M(l_2)$ valued random elements] to $\hat{u}$ and $x^\varepsilon \to x$, we have, as $\varepsilon \to 0$,

$$\hat{\mathcal{G}}^\epsilon\left(x^\varepsilon, \sqrt{\epsilon}\beta + \int_0^\cdot \hat{u}^\epsilon(s)\,ds\right) = \mathcal{G}^\epsilon(x^\varepsilon, \sqrt{\epsilon}B + \text{Int}(u_{\hat{u}^\epsilon}))$$
$$\to \mathcal{G}^0(x, \text{Int}(u_{\hat{u}}))$$
$$= \hat{\mathcal{G}}^0\left(x, \int_0^\cdot \hat{u}(s)\,ds\right).$$

Thus, part 2 of Assumption 2 is satisfied with $\mathcal{G}^\varepsilon$ replaced by $\hat{\mathcal{G}}^\varepsilon$, $\varepsilon \geq 0$. The result now follows on noting that if $\hat{I}_x(f)$ is defined by the right-hand side of (4.5) on replacing $\mathcal{G}^0$ there by $\hat{\mathcal{G}}^0$, then $\hat{I}_x(f) = I_x(f)$ for all $(x,f) \in \mathcal{E}_0 \times \mathcal{E}$. □

A. Budhiraja
V. Maroulas
Department of Statistics
and Operations Research
University of North Carolina
Chapel Hill, North Carolina 27599
USA
E-mail: budhiraj@email.unc.edu
       maroulas@email.unc.edu

P. Dupuis
Division of Applied Mathematics
Brown University
Providence, Rhode Island 02912
USA
E-mail: dupuis@dam.brown.edu